\documentclass[11pt]{article}

\usepackage[T1]{fontenc}
\usepackage{lmodern}
\usepackage[a4paper,margin=28mm]{geometry}
\usepackage{microtype}
\usepackage{amsmath,amssymb,amsthm,mathtools}
\usepackage{longtable}
\usepackage{csquotes}
\usepackage{xcolor}
\usepackage{colortbl}
\definecolor{colone}{rgb}{0.85,0.92,1}
\definecolor{refcolor}{rgb}{0,0,0.75}
\definecolor{citecolor}{rgb}{0,0.45,0.15}
\usepackage[colorlinks,linkcolor=refcolor,citecolor=citecolor,
 urlcolor=refcolor,pagebackref]{hyperref}
\renewcommand*{\backref}[1]{}
\renewcommand*{\backrefalt}[4]{%
 \ifcase #1 \or (p.~#2)\else (pp.~#2)\fi}

\allowdisplaybreaks

\newtheorem{theorem}{Theorem}[section]
\newtheorem{proposition}[theorem]{Proposition}
\newtheorem{lemma}[theorem]{Lemma}
\newtheorem{corollary}[theorem]{Corollary}
\newtheorem{observation}[theorem]{Observation}
\newtheorem{fact}[theorem]{Fact}
\newtheorem{conjecture}[theorem]{Conjecture}
\newtheorem{question}[theorem]{Question}
\theoremstyle{definition}
\newtheorem{definition}[theorem]{Definition}
\newtheorem{remark}[theorem]{Remark}
\theoremstyle{plain}
\DeclareMathOperator{\SL}{SL}
\DeclareMathOperator{\GL}{GL}
\DeclareMathOperator{\EL}{EL}
\DeclareMathOperator{\PSL}{PSL}
\DeclareMathOperator{\PGL}{PGL}
\DeclareMathOperator{\PSU}{PSU}
\DeclareMathOperator{\PSp}{PSp}
\DeclareMathOperator{\Sp}{Sp}
\DeclareMathOperator{\SU}{SU}
\DeclareMathOperator{\Aut}{Aut}
\DeclareMathOperator{\Out}{Out}
\DeclareMathOperator{\Alt}{Alt}
\DeclareMathOperator{\Sym}{Sym}
\DeclareMathOperator{\gap}{gap}
\DeclareMathOperator{\mindim}{mindim}
\DeclareMathOperator{\tr}{tr}
\DeclareMathOperator{\Cay}{Cay}
\DeclareMathOperator{\nrd}{nrd}
\DeclareMathOperator{\row}{row}
\DeclareMathOperator{\col}{col}
\DeclareMathOperator{\Reg}{Reg}
\DeclarePairedDelimiter{\angles}{\langle}{\rangle}
\DeclarePairedDelimiterXPP{\opnorm}[1]{}{\lVert}{\rVert}{_{\mathrm{op}}}{#1}

\newcommand{\F}{\mathbb F}

\newcommand{\Z}{\mathbb Z}

\newcommand{\R}{\mathbb R}

\newcommand{\C}{\mathbb C}

\newcommand{\POmega}{\mathrm P\Omega}

\newcommand{\defeq}{\coloneqq}

\title{Expansion in powers of simple groups of Lie type}
\author{Alexander Lubotzky\thanks{Weizmann Institute of Science, Rehovot.
\texttt{alex.lubotzky@mail.huji.ac.il},
\texttt{yotam.shomroni@weizmann.ac.il}.}
\and
Yotam Shomroni\footnotemark[1]}
\date{\today\\[6pt]\normalsize\textit{Dedicated to Avi Wigderson, an expansion leader, on his 70th birthday}}

\begin{document}

\maketitle

\begin{abstract}
We prove an expansion theorem for direct powers of all nonabelian finite
simple groups $G$ other than the alternating groups and the Suzuki groups
${}^2B_2(q)$.  If $G^m$ is
generated by $d\ge2$ elements, then $G^m$ has a symmetric generating multiset
of at most $Cd$ elements and normalized spectral gap at least
$C^{-1}d^{-\log_23}$,
for an absolute constant $C>0$.  In particular, this answers a question of Sabatini on
$\PSL_2(\F_p)^p$.  For $\PSL_2(\F_q)$, we use the expander graphs
constructed by Lubotzky--Samuels--Vishne, and an amplification trick for
direct powers; in high rank, we use the
Ershov--Jaikin-Zapirain property~(T) theorem.
\end{abstract}

\tableofcontents

\section{Introduction}

A basic problem on expansion of finite groups is the following.  Given a
family $\mathcal G$ of finite groups, are there $k$ and $\varepsilon>0$ such
that every $G\in\mathcal G$ has a multiset $A$ of at most $k$ elements with
$\Cay(G;A)$ an
$\varepsilon$-expander?\footnote{See the definition before
Theorem~\ref{thm:proved-global}.}  In
this case we say that $\mathcal G$ is an
\textbf{expanding family}.  Of course, a necessary condition is that the
minimal number of generators of $G$ is bounded over $\mathcal G$.
Various papers give positive and negative answers to this
question (see e.g.\ \cite{LubotzkyLie,KLN,BGT,BGGT,BK,Sabatini,Jez}).

One of the most famous examples is when $\mathcal G$ runs over all nonabelian
finite simple groups.
By \cite{SteinbergGen,AG} they are all $2$-generated, and by
\cite{KLN,BGT} they are also an expanding family, though it is not known
whether with $k=2$.\footnote{By \cite{BGGT}, for finite simple groups of Lie
type of bounded rank, the Cayley graph of a random pair of elements is an
expander with probability tending to $1$.}

In this paper we look at families whose members are direct powers $G^m$, the
direct product of $m$ copies of a nonabelian finite simple group $G$.  It is known that
very high powers of $G$ are still generated by few elements, see
\eqref{eq:hall-count}.

Sabatini~\cite[Question~4.5]{Sabatini} put forward the following test case:
is $\{\PSL_2(\F_p)^p\}_{p\text{ prime}}$ an expanding family?  In Corollary~\ref{cor:sabatini}
we answer this in the positive.
In fact, Theorem~\ref{thm:huge-powers} below gives even higher powers.  For
the rest of
the paper, let us denote by $\mathcal L$ the family of all nonabelian finite
simple groups excluding the alternating groups and the Suzuki groups
${}^2B_2(q)$.

\begin{theorem}
\label{thm:huge-powers}
The family $\mathcal G=\{G^{|G|}:G\in\mathcal L\}$ is an expanding family.
\end{theorem}

By the classification of finite simple groups, every group in $\mathcal L$
is a sporadic group or a group of Lie type other than ${}^2B_2(q)$.
For a finite group $G$, denote by $D_d(G)$ the set of its generating
$d$-tuples, and denote
\begin{equation}
 \label{eq:hall-count}
 \mu_d(G)\defeq\frac{|D_d(G)|}{|\Aut(G)|}.
\end{equation}
By a criterion of Hall~\cite{Hall} (see also \cite[Corollary~7]{KL}), the
largest $d$-generated power of a nonabelian finite simple group $G$ is
$G^{\mu_d(G)}$.  
As $|\Aut(G)|=|G|\cdot|\Out(G)|$, and since, by \cite{Dixon,KL,LiebeckShalev},
$|D_2(G)|/|G|^2\to 1$ as ${|G|\to\infty}$, we have
$\mu_2(G)\sim|G|/|\Out(G)|$ as $|G|\to\infty$.  In general,
$|G|^{d-2}\le\mu_d(G)\le|G|^{d-1}$ for $d\ge2$ (see
Observation~\ref{obs:multiplicity}).

We measure expansion by the normalized spectral gap, $\gap_G(S)$, of a finite
symmetric multiset $S$, defined in \eqref{eq:gap-def}: a multiset is
\textbf{symmetric} if every element and its inverse have the same multiplicity,
and \textbf{generating} if its support generates the group.\footnote{For a
symmetric multiset $S$ in a finite group $G$, the gap $\gap_G(S)$ is positive
if and only if $S$ is generating.}  All cardinalities
count multiplicity, and by $A\uplus B$ we denote the multiset in which the
multiplicities are summed.  We call $\Cay(G;A)$ an
\textbf{\boldmath$\varepsilon$-expander} when $\gap_G(A\uplus A^{-1})\ge\varepsilon$.

\begin{theorem}[Main theorem]
\label{thm:proved-global}
There is an absolute constant $C>0$ with the following
property.  Let $G\in\mathcal L$, let $d\ge2$ and $m\ge1$, and suppose that
$G^m$ is generated by $d$ elements.  Then $G^m$ has a symmetric generating
multiset $S$ such that
\[
 |S|\le Cd,
 \qquad
 \gap_{G^m}(S)\ge C^{-1}d^{-\log_23}.
\]
\end{theorem}

We stress that the size of the generating set and the spectral gap are
independent of the group $G$, the field, the characteristic, and the
rank.\footnote{For every $\delta>0$, the same proof gives
$\gap_{G^m}(S)\ge C^{-1}d^{-1-\delta}$ for all but finitely many
$G\in\mathcal L$; see Remark~\ref{rem:exponent}.}

In a companion paper~\cite[Theorem~1.1]{Suz}, Theorem~\ref{thm:proved-global}
is proved for the Suzuki groups ${}^2B_2(q)$.  The Suzuki groups do not have
the subgroups isomorphic to $\SL_2(\F_q)$ that our proof in bounded rank uses,
so our method does not suffice.

Expansion of powers of simple groups was studied by Kassabov.  By
\cite[Corollary~3.1]{KassabovSym}, for every $s>6$ the group
$\SL_{3s}(\F_2)^{(2^{3s}-1)^5}$ has a generating set of $26$
involutions whose Cayley graph is an $\varepsilon$-expander for an absolute
$\varepsilon$.
Lemma~\ref{lem:ring-generation} below, which for $q=2$ follows from the
proof of \cite[Lemma~2.1]{KassabovSym} (in its arXiv version), applies
to every power $m$, at
the cost of more generators.
In the following table, each row says that every $d$-generated power $G^m$
has a symmetric generating multiset $S$ with $|S|\le Cd$ and
$\gap_{G^m}(S)\ge C^{-1}d^{-\alpha}$, for the exponent $\alpha$ and the
constant $C$ of that row.
\begin{center}
\phantomsection\label{tab:known-results}%
\begin{tabular}{@{}l@{\ }|@{\,}c@{\,}|@{\ }l@{\ }|@{\ }l@{}}
 $G$ & $\alpha$ & $C$ & source\\ \hline
 $\SL_{3n}(\F_q)$ & $2$ & absolute
  & \cite{KassabovSym} with \cite{EJ}\\
 untwisted of rank $\ge2$, not $B_2,C_N,G_2$
  & $2$ & depends on the rank
  & \cite{EJK}\\
 every $G\in\mathcal L$ & $\log_23$ & absolute
  & Theorem~\ref{thm:proved-global}\\
\end{tabular}
\end{center}
The first two rows are not stated in the cited works: we show in
Corollary~\ref{cor:type-A-powers} and Remark~\ref{rem:fixed-type} how they
follow from the cited theorems.  The two main novelties in our work are
Theorem~\ref{thm:amplification} (amplification) and the construction of
expanders for high powers of $\PSL_2(\F_q)$.  This case
is harder in the sense that, by \cite[Theorem~9.7]{EJK}, a group with
property~(T) maps onto $\PSL_2(\F_q)$ for only finitely many $q$.
It is
possible, however, that Sabatini's question was meant to be the easiest case.

It is interesting to compare Theorem~\ref{thm:proved-global} with the
following recent result.

\begin{theorem}[{Salehi Golsefidy--Li~\cite[Theorem~1]{SGL}}]
\label{thm:SGL}
Let $n_0\ge3$ and let $p\ge5$ be a prime.  There is a $C=C(n_0,p)$ with the
following property.  Let $d\ge1$ and suppose that $\SL_{n_0}(\F_p)^m$ is
generated by $d$ elements.  Then $\SL_{n_0}(\F_p)^m$ has a symmetric
generating multiset $S$ such that
\[
 |S|\le Cd,
 \qquad
 \gap_{\SL_{n_0}(\F_p)^m}(S)\ge C^{-1}.
\]
\end{theorem}

Salehi Golsefidy and Li~\cite[Section~5]{SGL} also asked whether
Theorem~\ref{thm:SGL} holds for $n_0=2$.
Theorems~\ref{thm:proved-global} and~\ref{thm:SGL} are very different and
complementary.  Theorem~\ref{thm:proved-global} works for every
$G\in\mathcal L$, with $G=\PSL_2(\F_q)$ the most difficult case and one that
Theorem~\ref{thm:SGL} does not cover,
while Theorem~\ref{thm:SGL}
impressively gives a gap independent of the number of generators.  The
combination of the two suggests the following conjecture.

\begin{conjecture}
There is an absolute constant $C>0$ such that every
$d$-generated power $G^m$ of every nonabelian finite simple group $G$ has a
symmetric generating multiset $S$ with $|S|\le Cd$ and
$\gap_{G^m}(S)\ge C^{-1}$.
\end{conjecture}

For the alternating groups, which $\mathcal L$ excludes, Bartholdi and
Kassabov~\cite{BK} built, for each fixed $a$, a group with property~(T)
surjecting onto $\Alt(n)^{f(n)}$ for infinitely many $n$, where
$\log f(n)\approx(\log n)^a$, and asked how large $f$ may be.  The growth
$(\log n)^a$ of $\log f(n)$ is still far below
$\log \mu_2(\Alt(n))\sim n\log n$.

For every nonabelian finite simple group $G$, at the extreme power $m=\mu_d(G)$, 
a generating $d$-tuple $A$ of $G^m$ is unique up to
the action of $\Aut(G^m)$ (see Observation~\ref{obs:unique-tuple}), so, up
to isomorphism, $G^m$ has only one Cayley graph with respect to $d$
generators.

\begin{question}
\label{q:maximal-power}
Fix $N\ge1$ and $d\ge2$.  Is there an $\varepsilon=\varepsilon(N,d)>0$ such
that the unique Cayley graph of $G^{\mu_d(G)}$ with respect to $d$ generators is an
$\varepsilon$-expander for every finite simple group $G$ of Lie type of
absolute Lie rank $N$ (see Section~\ref{sec:preliminaries})?
\end{question}

If Question~\ref{q:maximal-power} has an affirmative answer, then
$\Cay(G;S)$ is an $\varepsilon$-expander for every such $G$ and every
generating $d$-element subset $S$.
Becker and Breuillard proved~\cite[Theorem~1.1]{BB} that, apart from a
density zero set of primes
$p$, every Cayley graph with respect to $d$ generators of every finite
simple group of Lie type of bounded rank over $\F_{p^f}$ with $f$ bounded is
an $\varepsilon$-expander, for an $\varepsilon>0$ depending only on $d$ and
those bounds.  In
Remark~\ref{rem:upper-bound} we show that the $\varepsilon(N,d)$ from
Question~\ref{q:maximal-power} must be less than $2/(d(N+1))$ for $N\ge3$.
There we also describe a stronger result of Jezernik~\cite[Theorem~1.13]{Jez}
for the alternating groups, and show that a question of
his~\cite[Question~1.12]{Jez} is equivalent to
the analogue of Theorem~\ref{thm:proved-global} for the alternating groups.

For a finite group $G$, denote by $\mindim(G)$
the minimal
dimension of a nontrivial irreducible representation of $G$.  The
proof of Theorem~\ref{thm:proved-global} uses the
following amplification theorem for direct powers, which may be of
independent interest.  For a group $G$ and $m,n\in\Z_{\ge1}$, index the
$mn$ coordinates of $G^{mn}$ by the pairs $(i,j)$ with $1\le i\le m$ and
$1\le j\le n$, so that they form an $m$-by-$n$ table, and 
define two embeddings, $\row\colon G^m\hookrightarrow G^{mn}$ 
and $\col\colon G^n\hookrightarrow G^{mn}$, 
by $\row(a_1,\ldots, a_m)_{i, j}=a_i$ 
and $\col(b_1,\ldots, b_n)_{i, j}=b_j$.
For multisets $A$ in $G^m$ and $B$ in $G^n$, denote by $\row(A)$ and
$\col(B)$ their images in $G^{mn}$, and denote
\[
 \Sigma_{A,B}\defeq\row(A)\uplus\col(B)\subseteq G^{mn}.
\]

\begin{theorem}[Amplification]
\label{thm:amplification}
Let $G$ be a nontrivial finite perfect group, and let $A$ and $B$ be
symmetric generating multisets for
$G^m$ and $G^n$, with
$|A|\le |B|$ and
gaps $\varepsilon_A,\varepsilon_B>0$.  Then
$\Sigma_{A,B}$ is a symmetric generating multiset of $G^{mn}$ of size
$|A|+|B|$, and
\begin{equation}
 \label{eq:amplification-gap}
 \gap_{G^{mn}}(\Sigma_{A,B})
 \ge\frac{|A|}{|A|+|B|}
 \left(1-\frac1{\mindim(G)}\right)\min\{\varepsilon_A,\varepsilon_B\}
 \ge\frac{|A|}{|B|}\,\frac{\min\{\varepsilon_A,\varepsilon_B\}}4,
\end{equation}
and the denominator $4$ may be replaced by $3$ when $G$ is a nonabelian
finite simple group.
\end{theorem}

\subsection{Proof overview}

\paragraph{Amplification.}
The proof of Theorem~\ref{thm:amplification} uses
Corollary~\ref{obs:angle-mindim}: in every nontrivial irreducible
representation of $G^{mn}$, the orthogonal projection onto the
$\col(G^n)$-fixed vectors, restricted to the $\row(G^m)$-fixed vectors, has
norm at most $1/\mindim(G)$.  In the terms defined by Ershov and
Jaikin-Zapirain~\cite[Section~2.1]{EJ}, the two subspaces of fixed vectors
are therefore $1/\mindim(G)$-orthogonal, and we deduce a spectral gap as
in~\cite[Lemma~3.1(b)]{EJ}.  Rozenman, Shalev and
Wigderson~\cite[Section~5.1]{RSW} studied the expansion of $G\times G$ with
respect to generators of the form $(y,y^{-1})$, whose two coordinates are
correlated, as are the coordinates of each element of $\Sigma_{A,B}$.
Observation~\ref{obs:row-column-generate} is similar
to~\cite[Observation~5.3]{RSW}.

With Theorem~\ref{thm:amplification} we increase the power: from
expanding sets (that is, symmetric multisets with positive gap) in
$G^m$ and $G^n$ we build one in $G^{mn}$.  Iterating this construction $j$
times with $n=m$ replaces $m$ by $m^{2^j}$, multiplies the number of
generators by $2^j$, and, for a nonabelian finite simple group $G$, multiplies
the lower bound on the gap by $3^{-j}$.  In this way we prove
Corollary~\ref{cor:reach}: if $G^{m_0}$, for some $m_0\ge|G|^{1/\beta}$ with
$\beta\ge1$, has an expanding set of size at most $s$ and gap at least
$\eta$, then
every $d$-generated power $G^m$ has an expanding set of size at most
$2\beta sd$ and gap at least $\eta\,(2\beta d)^{-\log_23}$.

\medskip
The proof of Theorem~\ref{thm:proved-global} treats high rank and bounded
rank separately.

\paragraph{High rank.}
The expansion is obtained from one group with a Kazhdan-type bound, mapped
onto many simple quotients at once.  For an integer $r\ge2$, let
$R=\Z\angles{x_1,\ldots,x_r}$ be the free associative unital ring over
$\Z$, and let $\EL_3(R)$ be the group generated by the elementary matrices
$I_3+aE_{ij}$ with $a\in R$ and $i\ne j$.  Ershov and
Jaikin-Zapirain~\cite{EJ} proved a quantitative Kazhdan-constant estimate
for $\EL_3(R)$.  By Corollary~\ref{cor:standard-gap}, which follows from
it, the $12(r+1)$ matrices $I_3\pm aE_{ij}$ with $a\in\{1,x_1,\ldots,x_r\}$
and $i\ne j$ have spectral gap at least $2^{-21}(r+1)^{-2}$ in every finite
quotient.

Denote by $M_n(\F_q)$ the algebra of $n$-by-$n$ matrices over the field
$\F_q$.  By Lemma~\ref{lem:ring-generation}, for every
$m\le|M_n(\F_q)|^{r-2}$ there is a surjective ring homomorphism
$R\twoheadrightarrow M_n(\F_q)^m$, and it induces a surjection
$\EL_3(R)\twoheadrightarrow\SL_{3n}(\F_q)^m$.  For $m=1$, Kassabov used
such surjections in~\cite{KassabovUniv} to construct expanders of unbounded
rank.  We use only a fixed $r$, and obtain the higher powers by
amplification.

\paragraph{Bounded rank.}
For $\PSL_2(\F_q)$ we use the Cartwright--Steger group, constructed
in~\cite{CS}.  By a theorem of Lubotzky, Samuels and Vishne~\cite{LSV}
(Theorem~\ref{thm:LSV}), its congruence quotients
have Ramanujan Cayley graphs; the power is then raised by repeated
amplification.  The case $m=1$ is \cite[Theorem~3.7]{LubotzkyLie}, and we
follow the proof there.
For a prime power $q$, there is a congruence quotient of the
Cartwright--Steger group isomorphic to
\[
 \PSL_2(\F_q)^{m_q},\qquad
 m_q\defeq\begin{cases}(q-3)/2,&q\text{ odd},\\q-2,&q\text{ even}.
 \end{cases}
\]
The Cayley graph of $\PSL_2(\F_q)^{m_q}$ with the $q+1$ Cartwright--Steger
generators is a Ramanujan graph, so all eigenvalues of its averaging operator
on the functions orthogonal to the constants are at most $2\sqrt q/(q+1)$.  These generators form at most two orbits under
conjugation by the intersection of $\PSL_2(\F_q)^{m_q}$ with a diagonally
embedded nonsplit torus of $\PGL_2(\F_q)^{m_q}$, defined in the proof of
Proposition~\ref{prop:bounded-rank-one}.  By \cite[Theorem~3.7]{LubotzkyLie}
(Theorem~\ref{thm:lubotzky} below), whose proof uses Selberg's $3/16$
theorem~\cite{Selberg}, the group $\PSL_2(\F_q)$ has a generating multiset
of $8$ elements with gap at least an absolute constant $\varepsilon_0>0$.
As in that proof, we take the diagonal copy of this multiset in
$\PSL_2(\F_q)^{m_q}$, and one generator from each orbit with its inverse.
By Lemma~\ref{lem:conjugate-transfer}, their union, a generating multiset
of $\PSL_2(\F_q)^{m_q}$ of size at most $12$, has gap at least an absolute
constant.  For $q\ge7$, we have $m_q\ge q^{1/3}>|\PSL_2(\F_q)|^{1/9}$, so we
use Corollary~\ref{cor:reach} with $\beta=9$: every $d$-generated power
$\PSL_2(\F_q)^m$ has a symmetric generating multiset of $O(d)$ elements with
gap $\Omega(d^{-\log_23})$.

\paragraph{From type $A$ to the other types.}\phantomsection\label{par:other-types}
In both branches (high rank and bounded rank), each group $G$ is a product
of a bounded number of subgroups, each of which is a quotient of some
$\SL_{3n}(\F_{q_0})$ in high rank and of $\SL_2(\F_q)$ in bounded rank.
By \cite[Corollary~2.2]{LubotzkyLie}, if a finite group is a
product of $c$ subgroups, each with a generating set of at most $k$ elements
whose Cayley graph is a $\delta$-expander, then the union of these sets is a
generating set whose Cayley graph is an $\varepsilon$-expander, for an
$\varepsilon>0$ depending only on $c$, $k$ and $\delta$.
Lemma~\ref{lem:product-transfer} is a quantitative form of this corollary.

For high rank groups, by Lemma~\ref{lem:large-A}, which follows from
\cite[Theorem~1 and Section~4, Lemma~2]{Nikolov}, every finite simple
classical group $G$ of absolute Lie rank at least $33$ is a product of at most $2^{13}$ conjugates
of a subgroup $H$ with $|H|\ge|G|^{1/3}$, which is a quotient of some
$\SL_{3n}(\F_{q_0})$.  The number $33$ is the least rank for which the
estimates in the proof of Lemma~\ref{lem:large-A} give $|H|\ge|G|^{1/3}$;
smaller ranks are treated as bounded rank.  The same factorization holds in
every direct power
$G^m$, with $H^m$ in place of $H$, so Lemma~\ref{lem:product-transfer}
applies to $G^m$.  Ershov, Jaikin-Zapirain and
Kassabov~\cite[Theorem~9.3]{EJK} constructed one group with property~(T)
mapping onto every finite simple group of Lie type of twisted rank at least
two.  In Remark~\ref{rem:fixed-type} we explain how to obtain expanders in high
powers of the untwisted groups of rank at least two, not of type $B_2$,
$C_N$ or $G_2$, from property~(T) of elementary Chevalley groups over polynomial
rings in many variables.  The constants obtained in that way depend on the
rank, whereas ours are absolute.  By
Remark~\ref{rem:general-r}, for any number $r\ge3$ of variables of the free
ring $R$, our high-rank construction gives expanders in every $G^m$ with
$m\le|G|^{(r-2)/27}$, and, for $r=27(d-1)+2$, these powers include every
$d$-generated power, with gap $\Omega(d^{-2})$.  Therefore, in high rank, amplification is
needed only to improve the gap from $\Omega(d^{-2})$ to
$\Omega(d^{-\log_23})$.  For this improvement we choose $r=29$, for which
$|G|^{(r-2)/27}=|G|$, and obtain every $d$-generated power from $G^{|G|}$
by Corollary~\ref{cor:reach}.

For bounded rank groups, Liebeck, Nikolov and Shalev~\cite{LNS} showed that
every finite simple group of Lie type over $\F_q$ other than ${}^2B_2(q)$ is
a product of subgroups isomorphic to $\SL_2(\F_q)$ or $\PSL_2(\F_q)$, whose
number is bounded in terms of the rank.  The sporadic groups and
the groups of Lie type of bounded rank over the smallest fields form a finite
family, and we treat them separately using amplification.  The Suzuki groups
${}^2B_2(q)$ are treated in~\cite{Suz}: for $m\le q^{1/100}$, the Cayley
graph of ${}^2B_2(q)^m$ with respect to a random pair of elements has
spectral gap bounded below uniformly in $q$ and $m$, with probability
tending to $1$ as $q\to\infty$, as Breuillard, Green and Tao~\cite{BGT}
proved for $m=1$; then Corollary~\ref{cor:reach} applies.

\subsection*{Acknowledgements and AI disclosure}
The authors were partially supported by the European Research Council (ERC)
under the European Union's Horizon 2020 research and innovation programme
(grant agreement No.~882751, TeStability).

The mathematical ideas of this paper are the authors' own.  The
authors acknowledge the use of Anthropic's Claude Opus for reviewing and
refining proofs, editing, finding references and filling in some technical
details.  The
constants in Theorem~\ref{thm:amplification} were improved by OpenAI's
ChatGPT; our earlier bounds were slightly weaker, and depended on Ore's
conjecture, proved in~\cite{LOST}.

\section{Preliminaries}
\label{sec:preliminaries}

Denote $\PSL_n(\F_q)\defeq\SL_n(\F_q)/Z(\SL_n(\F_q))$, and define
$\PGL_n(\F_q)$, $\PSU_n(\F_{q^2})$, $\PSp_{2n}(\F_q)$ and
$\POmega_{2n}^{\pm}(\F_q)$ similarly, by dividing by the center.  Here
$\SU_n(\F_{q^2})$ is the group of determinant-one isometries of a
nondegenerate Hermitian form on $\F_{q^2}^n$.  For a nondegenerate quadratic
form on $\F_q^k$, the group $\Omega_k(\F_q)$ is the derived subgroup of its
orthogonal group; for $\Omega^+_{2n}(\F_q)$ and $\Omega^-_{2n}(\F_q)$, the
largest subspaces of $\F_q^{2n}$ on which the form vanishes identically have
dimension $n$ and $n-1$, respectively.  We use the standard Dynkin-type labels
for the finite simple groups of Lie type, where the subscript $N$ is the
\textbf{absolute Lie rank}, that is, the rank of the underlying root system:
\[
\begin{array}{ll}
 A_N(q)=\PSL_{N+1}(\F_q), & {}^2A_N(q)=\PSU_{N+1}(\F_{q^2}),\\
 B_N(q)=\Omega_{2N+1}(\F_q), & C_N(q)=\PSp_{2N}(\F_q),\\
 D_N(q)=\POmega_{2N}^+(\F_q), & {}^2D_N(q)=\POmega_{2N}^-(\F_q).
\end{array}
\]

Let $G$ be a finite group, let $\pi\colon G\to\GL(V)$ be a unitary
representation on a finite-dimensional complex Hilbert space $V$, and let
$V^G$ be the subspace of $G$-fixed vectors.  We denote by
\[
 P_G\defeq\frac1{|G|}\sum_{g\in G}\pi(g)
\]
the orthogonal projection onto $V^G$.  
A finite symmetric multiset $S$ in $G$
induces the following averaging operator $M_S^\pi$, to which we also attach a quadratic form $E_S^\pi$:
\[
 M_S^\pi\defeq\frac1{|S|}\sum_{s\in S}\pi(s),
 \qquad
 E_S^\pi(v)\defeq\angles{(I-M_S^\pi)v,v}
 \quad(\text{for }v\in V).
\]
Since $S$ is symmetric,
\begin{equation}\label{eq:dirichlet}
 E_S^\pi(v)=\frac1{2|S|}\sum_{s\in S}\|\pi(s)v-v\|^2.
\end{equation}
We denote by $\Reg_G$ the left regular representation of $G$ on
$\ell^2(G)$, so $(\Reg_G(g)v)(x)=v(g^{-1}x)$, and denote by $\mathbf1$
the subspace of constant functions.  The normalized one-sided spectral gap
is
\begin{equation}
 \label{eq:gap-def}
 \gap_G(S)\defeq
 \inf_{0\ne v\perp\mathbf 1}
 \frac{E_S^{\Reg_G}(v)}{\|v\|^2}.
\end{equation}
If $Q$ is a quotient of $G$ and $\bar S$ is the image of $S$ in $Q$, a
multiset, then, as the regular representation of $Q$ is a subrepresentation
of the regular representation of $G$,
\begin{equation}
 \label{eq:quotient-gap}
 \gap_Q(\bar S)\ge\gap_G(S).
\end{equation}
If $g=g_1\cdots g_k$ in $G$ and $v\in V$, then
$\pi(g)v-v=\sum_{l=1}^k\pi(g_1\cdots g_{l-1})(\pi(g_l)v-v)$, so
\begin{equation}\label{eq:product-displacement}
 \|\pi(g)v-v\|\le\sum_{l=1}^k\|\pi(g_l)v-v\|.
\end{equation}

\begin{observation}
\label{obs:gap-transfer}
If $S$ is a symmetric multiset in a finite group $G$ and
$\gap_G(S)\ge\varepsilon>0$, then for every unitary representation
$\pi\colon G\to\GL(V)$, every $v\in V$ and every $g\in G$,
\[
 E_S^\pi(v)\ge\varepsilon\,\|v-P_Gv\|^2
 \quad\text{and}\quad
 \|\pi(g)v-v\|\le2\|v-P_Gv\|\le2\sqrt{E_S^\pi(v)/\varepsilon}.
\]
\end{observation}

\begin{proof}
As $P_Gv$ is fixed by $G$, by the formula \eqref{eq:dirichlet} for $E^\pi$
we have $E_S^\pi(v)=E_S^\pi(v-P_Gv)$, so for the first inequality we may
assume that $v\perp V^G$.  Then $v$ is an orthogonal sum of vectors in
nontrivial irreducible subrepresentations of $\pi$, and $M_S^\pi$ preserves
each of these subrepresentations.  Every nontrivial irreducible
representation of $G$ occurs in $\Reg_G$ on $\mathbf1^\perp$, so the first
inequality follows from the definition \eqref{eq:gap-def} of the gap.
Moreover,
$\|\pi(g)v-v\|=\|\pi(g)(v-P_Gv)-(v-P_Gv)\|\le2\|v-P_Gv\|$, and, by the first
inequality, $2\|v-P_Gv\|\le2\sqrt{E_S^\pi(v)/\varepsilon}$.
\end{proof}

\begin{observation}
\label{obs:multiplicity}
Let $G$ be a nonabelian finite simple group and let $d\ge2$.  Then
$|G|^{d-2}\le\mu_d(G)\le|G|^{d-1}$.
Indeed, by \cite{SteinbergGen,AG}, the set $D_2(G)$ is nonempty (i.e.\ $G$ is
$2$-generated), and
$\Aut(G)$ acts freely on it, so $|D_2(G)|\ge|\Aut(G)|$.  As
$D_2(G)\times G^{d-2}\subseteq D_d(G)\subseteq G^d$ and $|\Aut(G)|\ge|G|$,
\[
 |G|^{d-2}\le\frac{|D_2(G)|\,|G|^{d-2}}{|\Aut(G)|}
 \le\frac{|D_d(G)|}{|\Aut(G)|}=\mu_d(G)\le\frac{|G|^d}{|G|}=|G|^{d-1}.
\]
\end{observation}

\begin{lemma}[{Kantor--Lubotzky~\cite[Proposition~6]{KL}}]
\label{lem:independent-simple}
Let $G$ be a nonabelian finite simple group, let $\Gamma$ be a finitely
generated group, and
let $\theta_1,\ldots,\theta_m\colon\Gamma\twoheadrightarrow G$ be surjective
homomorphisms.  The product map
$\Gamma\to G^m$ is surjective if and only if no two of the maps differ by an
automorphism of $G$.
\end{lemma}

\begin{observation}
\label{obs:unique-tuple}
Let $G$ be a nonabelian finite simple group, let $d\ge2$ and denote
$\mu\defeq\mu_d(G)$.  For every generating $d$-tuple $(y_1,\ldots,y_d)$ of
$G^\mu$, each $\Aut(G)$-orbit in $D_d(G)$ contains exactly one of the tuples
$((y_1)_i,\ldots,(y_d)_i)$ for $i\le\mu$.  Hence $\Aut(G^\mu)$ acts
transitively on the generating $d$-tuples of $G^\mu$.
Indeed, $\Aut(G)$ acts freely on $D_d(G)$, so, by the definition
\eqref{eq:hall-count} of $\mu$, it has exactly $\mu$ orbits there.  By
Lemma~\ref{lem:independent-simple}, where $\Gamma$ is the free group on
$x_1,\ldots,x_d$ and $\theta_i(x_k)=(y_k)_i$, the $\mu$ tuples are in
distinct orbits, hence one in each.  As $\Aut(G)\wr\Sym(\mu)$ acts on
$G^\mu$ by automorphisms, the transitivity follows.
\end{observation}

\subsection*{Bounded products of conjugates}

\begin{observation}
\label{obs:direct-power-product}
Let $G$ be a group, let $c,m\in\Z_{\ge1}$, and let $G_1,\ldots,G_c$ be
subgroups of $G$.  If
$G=G_1\cdots G_c$, then
$G^m=G_1^m\cdots G_c^m$, because the product in $G^m$ is coordinatewise.
\end{observation}

Lemma~\ref{lem:product-transfer} is a quantitative form
of~\cite[Corollary~2.2]{LubotzkyLie}.

\begin{lemma}[Expansion of a bounded product]
\label{lem:product-transfer}
Let $s,c\in\Z_{\ge1}$ and $\eta>0$, let $G$ be a finite group, and suppose
$G=G_1\cdots G_c$ for subgroups $G_1,\ldots,G_c\le G$.  If $T_i$ is a
symmetric multiset in $G_i$ with $|T_i|=s$ and $\gap_{G_i}(T_i)\ge\eta$ for
each $i$, then
\begin{equation}
 \label{eq:transfer-union}
 \gap_G\left(\biguplus_{i=1}^cT_i\right)
 \ge\frac{\eta}{2c^2}.
\end{equation}
\end{lemma}

\begin{proof}
The idea is that a unit vector moved little by every $T_i$ is close to the
fixed space of every $G_i$; as $G=G_1\cdots G_c$, it is then moved little by
every element of $G$, so it is not orthogonal to the constants.

Denote $A\defeq\biguplus_{i=1}^cT_i$, a multiset of size $cs$, and
$\pi\defeq\Reg_G$, and let $v\in\mathbf1^\perp$ be a unit vector.
As $v\perp\mathbf1$, we have $\sum_{g\in G}\pi(g)v=0$, so
\[
 \frac1{|G|}\sum_{g\in G}\|\pi(g)v-v\|^2
 =2\|v\|^2-2\operatorname{Re}\angles*{
  \frac{1}{|G|}\sum_{g\in G} \pi(g) v
  ,v}=2,
\]
so some $g\in G$ has $\|\pi(g)v-v\|\ge\sqrt2$.  Then $g=h_1\cdots h_c$ for
some $h_i\in G_i$.  By \eqref{eq:product-displacement}, by
Observation~\ref{obs:gap-transfer}, applied to the restriction of $\pi$ to
each $G_i$, and by Cauchy--Schwarz,
\[
 \sqrt2\le\|\pi(g)v-v\|
 \le\sum_{i=1}^c\|\pi(h_i)v-v\|
 \le\frac2{\sqrt\eta}\sum_{i=1}^c\sqrt{E_{T_i}^\pi(v)}
 \le\frac2{\sqrt\eta}\sqrt{c\sum_{i=1}^cE_{T_i}^\pi(v)},
\]
so $\sum_iE_{T_i}^\pi(v)\ge\eta/(2c)$, and, by the formula
\eqref{eq:dirichlet} for $E^\pi$,
\[
 E_A^\pi(v)
 =\frac1{2cs}\sum_{i=1}^c\sum_{t\in T_i}\|\pi(t)v-v\|^2
 =\frac1c\sum_{i=1}^cE_{T_i}^\pi(v)
 \ge\frac{\eta}{2c^2}.
\]
As $v\in\mathbf1^\perp$ is an arbitrary unit vector,
\eqref{eq:transfer-union} follows from the definition \eqref{eq:gap-def} of
the gap.
\end{proof}

We use Lemma~\ref{lem:conjugate-transfer} in the proof of
Proposition~\ref{prop:bounded-rank-one}.  The argument of
Lemma~\ref{lem:conjugate-transfer} was used in \cite[proof of
Theorem~3.7]{LubotzkyLie}, for $\PSL_2(\F_q)$ and its subgroup
$\PSL_2(\F_p)$, where $q$ is a power of a prime $p$.

\begin{lemma}
\label{lem:conjugate-transfer}
Let $\varepsilon_K,\varepsilon_G>0$, let $K\le G$ be finite groups, and let
$S_K$ and $S_G$ be symmetric multisets in $K$ and in $G$, respectively, with
$\gap_K(S_K)\ge\varepsilon_K$ and $\gap_G(S_G)\ge\varepsilon_G$.  Let $S_G'$
be a symmetric multiset in $G$, and suppose that, for every $x\in S_G$, there
are $t\in K$ and $x'\in S_G'$ such that $x=t^{-1}x't$.  Then
$U\defeq S_K\uplus S_G'$ satisfies
$\gap_G(U)\ge\varepsilon_G/\bigl((16/(|S_K|\,\varepsilon_K)+2)\,|U|\bigr)$.
\end{lemma}

\begin{proof}
The idea is that a unit vector moved little by $U$ is close to the fixed space
of $K$, because $S_K\subseteq U$; hence it is moved little by every
$x=t^{-1}x't\in S_G$, because $x'\in U$, and, as
$\gap_G(S_G)\ge\varepsilon_G$, it is not orthogonal to the constants.

Denote $\pi\defeq\Reg_G$, and let $v\in\mathbf1^\perp$ be a unit vector.  As
$S_K\subseteq U$,
\begin{equation}
 \label{eq:union-energy}
 2|U|\,E_U^\pi(v)=\sum_{u\in U}\|\pi(u)v-v\|^2
 \ge\sum_{s\in S_K}\|\pi(s)v-v\|^2=2|S_K|\,E_{S_K}^\pi(v).
\end{equation}
Let $x\in S_G$, and let $t\in K$ and $x'\in S_G'$ be such that
$x=t^{-1}x't$.  By \eqref{eq:product-displacement}, and by
Observation~\ref{obs:gap-transfer} for the restriction of $\pi$ to $K$,
applied to $t^{-1}$ and to $t$,
\[
 \|\pi(x)v-v\|
 \le\|\pi(t^{-1})v-v\|+\|\pi(x')v-v\|+\|\pi(t)v-v\|
 \le2\cdot2\sqrt{E_{S_K}^\pi(v)/\varepsilon_K}+\|\pi(x')v-v\|.
\]
By \eqref{eq:union-energy}, and as $x'\in U$, we have
$|S_K|\,E_{S_K}^\pi(v)\le|U|\,E_U^\pi(v)$ and
$\|\pi(x')v-v\|^2\le2|U|\,E_U^\pi(v)$.  Hence, since $(a+b)^2\le2a^2+2b^2$ for
all $a,b\in\R$,
\begin{equation}
 \label{eq:conjugate-bound}
 \|\pi(x)v-v\|^2
 \le32E_{S_K}^\pi(v)/\varepsilon_K+2\|\pi(x')v-v\|^2
 \le\bigl(32/(|S_K|\,\varepsilon_K)+4\bigr)|U|\,E_U^\pi(v).
\end{equation}
By the definition \eqref{eq:gap-def} of the gap, and by
\eqref{eq:conjugate-bound}, summed over $x\in S_G$,
\[
 \varepsilon_G\le E_{S_G}^\pi(v)
 =\frac1{2|S_G|}\sum_{x\in S_G}\|\pi(x)v-v\|^2
 \le\bigl(16/(|S_K|\,\varepsilon_K)+2\bigr)|U|\,E_U^\pi(v).
\]
As $v\in\mathbf1^\perp$ is an arbitrary unit vector, the lemma follows from
the definition \eqref{eq:gap-def} of the gap.
\end{proof}

\section{Amplification}

Throughout this section $G$ is a nontrivial finite perfect group,
$m,n\in\Z_{\ge1}$, the $mn$ coordinates of $G^{mn}$ are indexed by the pairs
$(i,j)$ with $1\le i\le m$ and $1\le j\le n$, and we denote
$G_{\row}\defeq\row(G^m)\cong G^m$ and $G_{\col}\defeq\col(G^n)\cong G^n$
(recall the embeddings \enquote{$\row$} and \enquote{$\col$} defined before
Theorem~\ref{thm:amplification}).  We denote by
$G_{ij}\cong G$ the direct factor of $G^{mn}$ in coordinate $(i,j)$.  For
$x\in G$, we denote by $\col_j(x)\in G_{\col}$ the element that is $x$ in
every coordinate of column $j$ and the identity elsewhere, and we define
$\row_i(x)\in G_{\row}$ similarly (see Figure~\ref{fig:angle}).

We denote $[x,y]\defeq xyx^{-1}y^{-1}$.  The following observation is similar
to \cite[Observation~5.3]{RSW}, where Rozenman, Shalev and Wigderson showed that
for every $a,b\in G$ there are $z_1,z_2,z_3\in G$ for which
$([a,b],1)=\prod_{i=1}^3(z_i,z_i^{-1})$.

\begin{observation}
\label{obs:row-column-generate}
The subgroups $G_{\row}$ and $G_{\col}$ generate $G^{mn}$.  Indeed, $G$ is
perfect, and $[\row_i(x),\col_j(y)]$ is the element of $G^{mn}$ that is
$[x,y]$ in coordinate $(i,j)$ and the identity elsewhere.
\end{observation}

\begin{definition}[{Ershov--Jaikin-Zapirain~\cite[Section~2.1]{EJ}}]
\label{def:eps-orthogonal}
Two subspaces $U_1$ and $U_2$ of a Hilbert space are
\textbf{\boldmath$\varepsilon$-orthogonal} if
$|\angles{u_1,u_2}|\le\varepsilon\|u_1\|\,\|u_2\|$ for all $u_1\in U_1$ and
$u_2\in U_2$, that is, if $\opnorm{P_{U_2}|_{U_1}}\le\varepsilon$, where
$P_{U_2}$ is the orthogonal projection onto $U_2$.
\end{definition}

Note that
$\opnorm{P_{U_2}|_{U_1}}=\cos\sphericalangle(U_1,U_2)
\defeq\sup_{0\ne u_1\in U_1,\,0\ne u_2\in U_2}\cos\sphericalangle(u_1,u_2)$,
where $\sphericalangle(u_1,u_2)$ is the angle between $u_1$ and $u_2$, defined
by $\cos\sphericalangle(u_1,u_2)\defeq|\angles{u_1,u_2}|/(\|u_1\|\,\|u_2\|)$.
By \cite[Lemma~3.1(b)]{EJ}, if a group is
generated by two subgroups whose fixed subspaces are $\varepsilon$-orthogonal,
with $\varepsilon<1$, in every unitary representation without nonzero fixed
vectors, and $\delta$ bounds from below the Kazhdan
constants\footnote{See the definition before
Corollary~\ref{cor:standard-gap}.} of generating sets of the two subgroups, then their union has Kazhdan constant at least
$\delta\sqrt{(1-\varepsilon)/2}$.  We show in Lemma~\ref{lem:angle} and
Corollary~\ref{obs:angle-mindim} that the subgroups $G_{\row}$ and
$G_{\col}$ satisfy this with $\varepsilon=1/\mindim(G)$.

\begin{lemma}
\label{lem:angle}
Let $\pi\colon G^{mn}\to\GL(V)$ be an irreducible unitary representation.
Denote
$V=\bigotimes_{i,j}V_{ij}$ where $V_{ij}$ is an irreducible representation of
$G_{ij}$.  Then the orthogonal projection
$P_{G_{\col}}$ onto
the space $V^{G_{\col}}$ of $G_{\col}$-fixed vectors, restricted to the
space $V^{G_{\row}}$ of $G_{\row}$-fixed vectors, satisfies
\begin{equation}\label{eq:angle}
 \opnorm{P_{G_{\col}}|_{V^{G_{\row}}}}
 \le\frac1{\max_{i,j}\dim V_{ij}}.
\end{equation}
\end{lemma}

\begin{proof}

Let
$\pi_{ij}\colon G\to\GL(V_{ij})$ denote the irreducible representation on
$V_{ij}$, and fix a pair $(i_0,j_0)$ with $\dim V_{i_0j_0}$ maximal.  As
permuting columns is an automorphism of $G^{mn}$ that maps $G_{\row}$ and
$G_{\col}$ onto themselves, we may assume $j_0=1$.  For
$\vec a=(a_1,\ldots,a_m)\in G^m$, denote
\[
 V^{\col}_1\defeq\bigotimes_{i=1}^mV_{i,1},
 \qquad
 \pi^{\col}_1(\vec a)\defeq\bigotimes_{i=1}^m\pi_{i,1}(a_i),
 \qquad
 P^{\col}_1\defeq\frac1{|G|}\sum_{g\in G}\pi^{\col}_1(g,\ldots,g).
\]
Here $\pi^{\col}_1$ is an irreducible representation of $G^m$ on the tensor
factor $V^{\col}_1$ of $V$ indexed by the first column,
and $P^{\col}_1$ is the orthogonal projection onto the vectors of
$V^{\col}_1$ fixed by every $\pi^{\col}_1(g,\ldots,g)$.  Order the tensor
factors of $V$ with $V_{1,1},\ldots,V_{m,1}$ first, and denote by $\otimes\,I$
tensoring with the identity on the other factors.  In the example of
Figure~\ref{fig:angle}, where $m=2$ and $n=3$, we have
$V^{\col}_1=V_{1,1}\otimes V_{2,1}$.  For
$\vec a\in G^m$ and $g\in G$, as $\row(\vec a)_{ij}=a_i$ for every $j$,
\begin{align}
 \label{eq:row-on-column}
 \pi(\row(\vec a))
 &=\pi^{\col}_1(\vec a)\otimes\bigotimes_{j=2}^n\bigotimes_{i=1}^m\pi_{ij}(a_i),
 \\
 \label{eq:col-on-column}
 \pi(\col_1(g))&=\pi^{\col}_1(g,\ldots,g)\otimes I .
\end{align}
Hence, by \eqref{eq:col-on-column},
$P_{\col_1(G)}=\frac1{|G|}\sum_{g\in G}\pi(\col_1(g))=P^{\col}_1\otimes I$.

\begin{figure}[h]
\centering
\small
\renewcommand{\arraystretch}{1.25}
\newcommand{\sixcells}[6]{\begin{array}{|>{\columncolor{colone}}c|c|c|}
 \hline #1&#2&#3\\ \hline #4&#5&#6\\ \hline\end{array}}
\begin{tabular}{@{}c@{\qquad}c@{\qquad}c@{}}
 $\sixcells{a_1}{a_1}{a_1}{a_2}{a_2}{a_2}$
 &$\sixcells{b_1}{b_2}{b_3}{b_1}{b_2}{b_3}$
 &$\sixcells{g}{1}{1}{g}{1}{1}$\\[.6ex]
 (a) $\row(a_1,a_2)$&(b) $\col(b_1,b_2,b_3)$&(c) $\col_1(g)$
\end{tabular}
\caption{Example: $m=2$, $n=3$.}
\label{fig:angle}
\end{figure}

\emph{Step 1.}
We claim that, as operators on $V^{\col}_1$,
\begin{equation}
 \label{eq:row-average}
 \frac1{|G|^m}\sum_{\vec a\in G^m}
  \pi^{\col}_1(\vec a)^{-1}P^{\col}_1\pi^{\col}_1(\vec a)
 =\frac{\tr P^{\col}_1}{\dim V^{\col}_1}\,I .
\end{equation}
Indeed, the average on the left side of \eqref{eq:row-average} commutes
with $\pi^{\col}_1(\vec b)$ for every $\vec b\in G^m$, so, by Schur's lemma, it
is a scalar operator.  As conjugation preserves the trace, the
scalar is $\tr P^{\col}_1/\dim V^{\col}_1$.

Let $v\in V^{G_{\row}}$ with $\|v\|=1$ and let $\vec a\in G^m$.  By
\eqref{eq:row-on-column}, the factors of $\pi(\row(\vec a))$ outside the first
column cancel in the following equation:
\[
 \pi(\row(\vec a))^{-1}(P^{\col}_1\otimes I)\pi(\row(\vec a))
 =\bigl(\pi^{\col}_1(\vec a)^{-1}P^{\col}_1\pi^{\col}_1(\vec a)\bigr)
  \otimes I ,
\]
and $\pi(\row(\vec a))$ is unitary and fixes $v$, so
\[
 \angles[\big]{\bigl(\bigl(\pi^{\col}_1(\vec a)^{-1}P^{\col}_1
  \pi^{\col}_1(\vec a)\bigr)\otimes I\bigr)v,v}
 =\angles{(P^{\col}_1\otimes I)\pi(\row(\vec a))v,\pi(\row(\vec a))v}
 =\angles{(P^{\col}_1\otimes I)v,v}.
\]
Averaging over $\vec a\in G^m$ and using \eqref{eq:row-average}, and as
$P^{\col}_1\otimes I$ is an orthogonal projection,
\begin{equation}
 \label{eq:alpha}
 \|(P^{\col}_1\otimes I)v\|^2
 =\angles{(P^{\col}_1\otimes I)v,v}
 =\frac{\tr P^{\col}_1}{\dim V^{\col}_1}\,\|v\|^2
 =\frac{\tr P^{\col}_1}{\dim V^{\col}_1} .
\end{equation}

\emph{Step 2: the bound
$\tr P^{\col}_1/\dim V^{\col}_1\le(\dim V_{i_0,1})^{-2}$.}

Let $\chi_{=i_0}$ be the character of $\pi_{i_0,1}$ and $\chi_{\ne i_0}$ the
character of the representation
$g\mapsto\bigotimes_{i\ne i_0}\pi_{i,1}(g)$ of $G$, so that, for $g\in G$,
\[
 \chi_{=i_0}(1)=\dim V_{i_0,1},
 \qquad
 \dim V^{\col}_1=\chi_{=i_0}(1)\chi_{\ne i_0}(1),
 \qquad
 \tr\pi^{\col}_1(g,\ldots,g)=\chi_{=i_0}(g)\chi_{\ne i_0}(g).
\]
Hence, by the definition of $P^{\col}_1$,
\[
 \tr P^{\col}_1
 =\frac1{|G|}\sum_{g\in G}\tr\pi^{\col}_1(g,\ldots,g)
 =\frac1{|G|}\sum_{g\in G}\chi_{=i_0}(g)\chi_{\ne i_0}(g),
\]
which is the multiplicity in $\chi_{\ne i_0}$ of the complex conjugate
character $\bar\chi_{=i_0}$, when regarding the character $\chi_{\ne i_0}$ as a
sum of irreducible characters.  As $\bar\chi_{=i_0}$ is irreducible of degree
$\chi_{=i_0}(1)$, we have $\chi_{=i_0}(1)\tr P^{\col}_1\le\chi_{\ne i_0}(1)$,
so
\[
 \frac{\tr P^{\col}_1}{\dim V^{\col}_1}
 \le\frac{\chi_{\ne i_0}(1)}{\chi_{=i_0}(1)\dim V^{\col}_1}
 =\frac1{\chi_{=i_0}(1)^2}
 =\frac1{(\dim V_{i_0,1})^2}.
\]
Hence, by \eqref{eq:alpha},
$\|(P^{\col}_1\otimes I)v\|\le1/\dim V_{i_0,1}$.

As $\col_1(G)\le G_{\col}$ we have
$V^{G_{\col}}\subseteq V^{\col_1(G)}$, and projecting a vector onto a
larger subspace cannot shorten its image, so
$\|P_{G_{\col}}v\|\le\|P_{\col_1(G)}v\|\le1/\dim V_{i_0,1}$.  As $v$ was
an arbitrary unit vector of $V^{G_{\row}}$ and
$\dim V_{i_0,1}=\max_{i,j}\dim V_{ij}$, this proves \eqref{eq:angle}.
\end{proof}

\begin{corollary}
\label{obs:angle-mindim}
Let $\pi\colon G^{mn}\to\GL(V)$ be a unitary representation with no nonzero
$G^{mn}$-fixed vector.  Then
$\opnorm{P_{G_{\col}}|_{V^{G_{\row}}}}\le1/\mindim(G)$.
Indeed, $\pi$ decomposes orthogonally into nontrivial irreducible
representations on some spaces $W_1,\ldots,W_r$, each preserved by
$P_{G_{\row}}$ and $P_{G_{\col}}$, so $V^{G_{\row}}$ is the orthogonal sum of
the spaces $V^{G_{\row}}\cap W_k$.  Hence we may assume that $\pi$ is
irreducible and nontrivial.  Then, writing $V=\bigotimes_{i,j}V_{ij}$ with
$V_{ij}$ an irreducible representation of $G_{ij}$, as in
Lemma~\ref{lem:angle}, some $G_{ij}$ acts nontrivially on $V_{ij}$, so
$\max_{i,j}\dim V_{ij}\ge\mindim(G)$.
\end{corollary}

\begin{fact}
\label{obs:mindim}
Recall that $G$ is a nontrivial finite perfect group.  We have
$\mindim(G)\ge2$, because a representation of dimension $1$ has
abelian image and $G=[G,G]$.  If $G$ is a nonabelian finite simple group,
then $\mindim(G)\ge3$ (see \cite[Problem~3.3]{Isaacs}).  Indeed, suppose that $\tau\colon G\to\GL_2(\C)$ is a
nontrivial irreducible representation.  Its kernel is trivial and, as $G$ is
perfect, $\det\circ\,\tau$ is trivial, so $\tau$ embeds $G$ in $\SL_2(\C)$.
By a theorem of Frobenius, $2=\dim\tau$ divides $|G|$, so the group $G$ has
an element $h$ of order $2$.
Then $\tau(h)$ is an element of order $2$ in $\SL_2(\C)$, so $\tau(h)=-I$,
which is central; hence $h\in Z(G)=1$, a contradiction.
\end{fact}

Lemma~\ref{lem:far-subgroups} is the case of two subgroups of
\cite[Lemma~3.1(b)]{EJ}, with spectral gaps in place of Kazhdan constants in
order to keep the factor $|A|/(|A|+|B|)$, and its proof is the deduction
presented there.  For this lemma the assumption that $G$ is perfect is not
needed.

\begin{lemma}
\label{lem:far-subgroups}
Let $G$ be a finite group, let $H_1,H_2\le G$, and let
$\varepsilon\ge0$.  Suppose that, for every unitary representation
$\pi\colon G\to\GL(V)$ with no nonzero $G$-fixed vector, the
subspaces $V^{H_1}$ and $V^{H_2}$ are $\varepsilon$-orthogonal.  Let
$\varepsilon_A,\varepsilon_B>0$, and let $A$ and $B$ be symmetric multisets in
$H_1$ and in $H_2$, respectively, with $|A|\le|B|$,
$\gap_{H_1}(A)\ge\varepsilon_A$ and $\gap_{H_2}(B)\ge\varepsilon_B$.  Then
\[
 \gap_G(A\uplus B)
 \ge\frac{|A|\min\{\varepsilon_A,\varepsilon_B\}}{|A|+|B|}\,(1-\varepsilon).
\]
\end{lemma}

\begin{proof}
Denote $\pi\defeq\Reg_G$, let $v\in\mathbf1^\perp$ be a unit vector, where
$\mathbf1$ is the subspace of constant functions, and denote $x\defeq P_{H_1}v$
and $y\defeq P_{H_2}v$.  Then $x$ and $y$ are
vectors in $\mathbf1^\perp$ fixed by $H_1$ and by $H_2$, respectively, so, by
the assumption for the restriction of $\pi$ to $\mathbf1^\perp$,
\[
 |\angles{x,y}|
 \le\varepsilon\|x\|\,\|y\|
 \le\varepsilon\bigl(\|x\|^2+\|y\|^2\bigr)/2.
\]
As $\|x\|^2+\|y\|^2=\angles{x,v}+\angles{y,v}\le\|x+y\|\cdot\|v\|=\|x+y\|$,
it follows that
\[
 \bigl(\|x\|^2+\|y\|^2\bigr)^2\le\|x+y\|^2
 =\|x\|^2+\|y\|^2+2\operatorname{Re}\angles{x,y}
 \le(1+\varepsilon)\bigl(\|x\|^2+\|y\|^2\bigr),
\]
so $\|x\|^2+\|y\|^2\le1+\varepsilon$.  As
$\|v-x\|^2=\|v\|^2-2\operatorname{Re}\angles{v,x}+\|x\|^2=1-\|x\|^2$ and
$\|v-y\|^2=1-\|y\|^2$,
\begin{equation}\label{eq:projection-sum}
 \|v-x\|^2+\|v-y\|^2\ge1-\varepsilon.
\end{equation}
By Observation~\ref{obs:gap-transfer}, applied to the restrictions of $\pi$
to $H_1$ and to $H_2$,
\[
 E_A^\pi(v)\ge\varepsilon_A\|v-x\|^2,
 \qquad
 E_B^\pi(v)\ge\varepsilon_B\|v-y\|^2.
\]
As $|A|\le|B|$, by the formula \eqref{eq:dirichlet} for $E^\pi$ and by
\eqref{eq:projection-sum},
\begin{align*}
 E_{A\uplus B}^\pi(v)
 &=\frac{|A|E_A^\pi(v)+|B|E_B^\pi(v)}{|A|+|B|}
 \ge\frac{|A|\varepsilon_A\|v-x\|^2+|B|\varepsilon_B\|v-y\|^2}{|A|+|B|}\\
 &\ge\frac{|A|\min\{\varepsilon_A,\varepsilon_B\}}{|A|+|B|}
  \bigl(\|v-x\|^2+\|v-y\|^2\bigr)
 \ge\frac{|A|\min\{\varepsilon_A,\varepsilon_B\}}{|A|+|B|}\,(1-\varepsilon).
\end{align*}
As $v\in\mathbf1^\perp$ is an arbitrary unit vector, the lemma follows from
the definition \eqref{eq:gap-def} of the gap.
\end{proof}

\begin{proof}[Proof of Theorem~\ref{thm:amplification}]
The multiset $\Sigma_{A,B}=\row(A)\uplus\col(B)$ is symmetric of size
$|A|+|B|$, and by Observation~\ref{obs:row-column-generate} it is
generating.  As $\row$ and $\col$ are isomorphisms onto $G_{\row}$ and
$G_{\col}$, the gaps of $\row(A)$ in $G_{\row}$ and of $\col(B)$ in
$G_{\col}$ are $\varepsilon_A$ and $\varepsilon_B$.  By
Corollary~\ref{obs:angle-mindim}, for every unitary representation
$\pi\colon G^{mn}\to\GL(V)$ with no nonzero $G^{mn}$-fixed vector, the
subspaces $V^{G_{\row}}$ and $V^{G_{\col}}$ are $1/\mindim(G)$-orthogonal.
Hence, by Lemma~\ref{lem:far-subgroups} for the group $G^{mn}$ and the
multisets $\row(A)$ and $\col(B)$, the first inequality of
\eqref{eq:amplification-gap} holds.  For the second inequality,
$|A|+|B|\le2|B|$, and, by Fact~\ref{obs:mindim}, $1-1/\mindim(G)\ge1/2$, and
$1-1/\mindim(G)\ge2/3$ if $G$ is nonabelian simple.
\end{proof}

Let us now fix, for the following
Corollaries~\ref{cor:square-simple}--\ref{cor:reach}, a nonabelian finite
simple group $G$, integers $m_0,s\ge1$, and a symmetric generating multiset
$A$ of $G^{m_0}$ of size at most $s$ and gap at least $\eta>0$.

\begin{corollary}
\label{cor:square-simple}
The multiset $\Sigma_{A,A}$ is a symmetric generating multiset of
$G^{m_0^2}$ of size at most $2s$ and gap at least $\eta/3$.  Indeed, take
$A=B$ in Theorem~\ref{thm:amplification}.
\end{corollary}

\begin{corollary}
\label{obs:balanced-iteration}
For every integer $j\ge0$ and every integer $m$ with $1\le m\le m_0^{2^j}$,
the group $G^m$ has a symmetric generating multiset of size at most $2^js$
and gap at least $3^{-j}\eta$.  Indeed, by Corollary~\ref{cor:square-simple},
applied $j$ times, the group $G^{m_0^{2^j}}$ has such a multiset, and, by
\eqref{eq:quotient-gap}, so does its quotient $G^m$.
\end{corollary}

\begin{corollary}
\label{cor:reach}
Let $\beta\in\R_{\ge1}$, and suppose that $m_0\ge|G|^{1/\beta}$.  Then every
$d$-generated power $G^m$ has a symmetric generating multiset $S$ with
\[
 |S|\le2\beta s\,d,
 \qquad
 \gap_{G^m}(S)\ge\eta\,(2\beta d)^{-\log_23}.
\]
\end{corollary}

\begin{proof}
As $G^m$ is nonabelian, $d\ge2$, so $\beta(d-1)\ge1$.  Choose
$j\in\Z_{\ge0}$ such that $\beta(d-1)\le2^j<2\beta d$.  By
Observation~\ref{obs:multiplicity} and since $m_0\ge|G|^{1/\beta}$,
\[
 m\le\mu_d(G)\le|G|^{d-1}\le m_0^{\beta(d-1)}\le m_0^{2^j}.
\]
By Corollary~\ref{obs:balanced-iteration}, the group $G^m$ has a
symmetric generating multiset of size at most $2^js\le2\beta s\,d$ and gap at
least $3^{-j}\eta$, and $3^{-j}=(2^j)^{-\log_23}>(2\beta d)^{-\log_23}$.
\end{proof}

\section{High rank}
\label{sec:high-rank}

We construct expanders in powers of $\SL_{3n}(\F_q)$ first
(Proposition~\ref{prop:type-A}), and pass, by Theorem~\ref{thm:nikolov}, to
the other classical groups of absolute Lie rank at least $33$
(Proposition~\ref{prop:high-rank-range}).

\subsection*{Noncommutative universal lattices and
\texorpdfstring{type-$A$}{type-A} quotients}

For an integer $r\ge2$, denote by $R\defeq\Z\angles{x_1,\ldots,x_r}$ the
free associative unital ring.  For a unital ring $A$ and an integer $k\ge2$,
denote
\[
 e_{ij}(a)\defeq I_k+aE_{ij}
 \quad\text{for $a\in A$ and distinct $i,j\in\{1,\ldots,k\}$},
\]
where $E_{ij}$ is the matrix with a single $1$ in position $(i,j)$, and $0$
elsewhere, and denote by $\EL_k(A)$ the group generated by these matrices.
Following \cite{EJ}, we
call $\EL_3(R)$ a \textbf{noncommutative universal lattice}.  By
\cite[equation~(6.1)]{EJ}, it is generated by the symmetric set
$Y\defeq\{e_{ij}(a)^{\pm1}:i\ne j,\ a\in\{1,x_1,\ldots,x_r\}\}$, of size
$|Y|=12(r+1)$.

A finite subset $X$ of a group $\Gamma$ has \textbf{Kazhdan constant} at least
$\kappa$ if, for every unitary representation $\pi\colon\Gamma\to\GL(V)$
without nonzero fixed vectors, every unit vector $v\in V$ satisfies
$\max_{x\in X}\|\pi(x)v-v\|\ge\kappa$.  By \cite[equation~(6.2)]{EJ} for
$\EL_3(R)$, the subset $\{e_{ij}(a):i\ne j,\ a\in\{1,x_1,\ldots,x_r\}\}$ of
$Y$, and hence $Y$ itself, has Kazhdan constant at least
\[
 \kappa_r\defeq\frac1{8(12\sqrt{2r}+6+36\sqrt2)}.
\]

\begin{corollary}
\label{cor:standard-gap}
Let $Q$ be a finite quotient of $\EL_3(R)$, and let the multiset $T$ be the
image of $Y$ in $Q$.  Then $\gap_Q(T)\ge2^{-21}(r+1)^{-2}$.
Indeed, the restriction of $\Reg_Q$ to $\mathbf1^\perp$ is a unitary
representation of $\EL_3(R)$ without nonzero fixed vectors, so for every
unit vector $v\perp\mathbf1$ some $y\in Y$ satisfies
$\|\Reg_Q(y)v-v\|\ge\kappa_r$.  As $Y$ is symmetric and
$\|\Reg_Q(y^{-1})v-v\|=\|v-\Reg_Q(y)v\|$, by the formula
\eqref{eq:dirichlet} for $E_T^{\Reg_Q}$, and as $r\ge2$,
\[
 E_T^{\Reg_Q}(v)
 =\frac1{2|T|}\sum_{t\in T}\|\Reg_Q(t)v-v\|^2
 \ge\frac{\kappa_r^2}{|Y|}
 =\frac{1}{12(r+1)\cdot64\,(12\sqrt{2r}+6+36\sqrt2)^2}
 \ge\frac1{2^{21}(r+1)^2}.
\]
\end{corollary}

The following lemma is a special case of the analogue of
Lemma~\ref{lem:independent-simple} for rings.

\begin{lemma}[{Kravchenko--Mazur--Petrenko~\cite[Theorem~6.1]{KMP}}]
\label{lem:independent-rings}
Let $q$ be a prime power, let $n,k\ge1$, and let
$\theta_1,\ldots,\theta_m\colon\Z\angles{x_1,\ldots,x_k}\twoheadrightarrow
M_n(\F_q)$ be surjective ring homomorphisms.  The product map
$\Z\angles{x_1,\ldots,x_k}\to M_n(\F_q)^m$ is surjective if and only if no
two of the maps differ by an automorphism of the ring $M_n(\F_q)$.
\end{lemma}

In \cite{KMP}, the lemma is stated for the $m$th power of any
finite-dimensional simple algebra over a field $F$, and for generation as an
$F$-algebra.  The case above is $F=\F_p$, where $p$ is the characteristic of
$\F_q$: every subring of $M_n(\F_q)^m$ is an $\F_p$-subalgebra, and every ring
automorphism of $M_n(\F_q)$ is $\F_p$-linear.

For $q=2$, the next lemma follows from the proof of
\cite[Lemma~2.1]{KassabovSym} (in its arXiv version); for every $q$ it follows from
Lemma~\ref{lem:independent-rings} and \cite[Theorem~6.4]{KMP}.

\begin{lemma}
\label{lem:ring-generation}
Let $q$ be a prime power and let $n,t\ge1$.  Then $M_n(\F_q)^{q^{n^2t}}$ is a
quotient of the ring $\Z\angles{x_1,\ldots,x_{t+2}}$.
\end{lemma}

\begin{proof}
Let $p$ be the characteristic of $\F_q$.  By \cite[Theorem~6.4]{KMP}, a
simple separable algebra over a field $F$ can be generated by two elements as
an $F$-algebra.  The $\F_p$-algebra $M_n(\F_q)$ is simple, and it is
separable as $\F_p$ is perfect, so there are $\alpha,\beta\in M_n(\F_q)$ that
generate $M_n(\F_q)$ as an $\F_p$-algebra, that is, as a ring.  For $w=(w_1,\ldots,w_t)\in M_n(\F_q)^t$, let
$\theta_w\colon\Z\angles{x_1,\ldots,x_{t+2}}\to M_n(\F_q)$ be the ring
homomorphism sending $x_1$ to $\alpha$, $x_2$ to $\beta$, and $x_{2+i}$ to $w_i$ for
$1\le i\le t$; it is surjective because $\alpha$ and $\beta$ are in its
image.  If $w,w'\in M_n(\F_q)^t$ and an automorphism $\sigma$ of $M_n(\F_q)$
satisfies
$\sigma\circ\theta_{w'}=\theta_w$, then $\sigma$ fixes $\alpha$ and $\beta$,
so $\sigma$ is the identity and $w'=w$.  Hence, by
Lemma~\ref{lem:independent-rings}, the product of the $q^{n^2t}$ maps
$\theta_w$ is a surjection
$\Z\angles{x_1,\ldots,x_{t+2}}\twoheadrightarrow M_n(\F_q)^{q^{n^2t}}$.
\end{proof}

\begin{proposition}[A type-$A$ quotient]
\label{prop:type-A}
Let $q$ be a prime power, let $n\ge1$, and assume $r\ge3$.  If
$1\le m\le|\SL_{3n}(\F_q)|^{(r-2)/9}$, then there is a surjection
$\theta\colon\EL_3(R)\twoheadrightarrow\SL_{3n}(\F_q)^m$, and the multiset
$\theta(Y)$ satisfies
\[
 |\theta(Y)|=12(r+1),
 \qquad
 \gap_{\SL_{3n}(\F_q)^m}(\theta(Y))
 \ge2^{-21}(r+1)^{-2}.
\]
\end{proposition}

\begin{proof}
As $|\SL_{3n}(\F_q)|<q^{9n^2}$, we have $m\le q^{n^2(r-2)}$, so, by
Lemma~\ref{lem:ring-generation} with $t=r-2$, followed by a projection onto
$m$ coordinates, there is a surjective ring homomorphism
$R\twoheadrightarrow M_n(\F_q)^m$.  As
$e_{ij}(a_1,\ldots,a_m)=(e_{ij}(a_1),\ldots,e_{ij}(a_m))$ and
$e_{ij}(0)=I_3$, as $\EL_3(M_n(\F_q))\cong\EL_{3n}(\F_q)$ by
\cite[Section~6.1, footnote~5]{EJ}, and as $\F_q$ is a field, this
homomorphism induces a surjection
\[
 \theta\colon\EL_3(R)\twoheadrightarrow\EL_3\bigl(M_n(\F_q)^m\bigr)
 =\EL_3\bigl(M_n(\F_q)\bigr)^m\cong\EL_{3n}(\F_q)^m=\SL_{3n}(\F_q)^m.
\]
Finally, $|\theta(Y)|=|Y|=12(r+1)$, and, by Corollary~\ref{cor:standard-gap},
$\theta(Y)$ has the stated gap.
\end{proof}

\begin{corollary}
\label{cor:type-A-powers}
Let $n\ge1$, let $q$ be a prime
power and let $d\ge2$.  For every $d$-generated power $\SL_{3n}(\F_q)^m$,
there is a symmetric generating multiset $S$ with $|S|\le2^7d$ and
$\gap_{\SL_{3n}(\F_q)^m}(S)\ge2^{-28}d^{-2}$.
Indeed, $\PSL_{3n}(\F_q)^m$ is $d$-generated too, so
$m\le\mu_d(\PSL_{3n}(\F_q))\le|\PSL_{3n}(\F_q)|^{d-1}\le|\SL_{3n}(\F_q)|^{d-1}$ by
Observation~\ref{obs:multiplicity}.  Hence, for $r=9(d-1)+2$, the multiset
$\theta(Y)$ in Proposition~\ref{prop:type-A} is symmetric, of size
$12(9d-6)\le2^7d$ and gap at least $2^{-21}(9d-6)^{-2}\ge2^{-28}d^{-2}$.
\end{corollary}

\subsection*{From type \texorpdfstring{$A$}{A} to the other classical types}

Let $G$ be a finite simple classical group over $\F_q$ of absolute Lie rank
$N$.  Then $G=\tilde G/Z(\tilde G)$ for some $\tilde G$ in the following table;
let $q_0$ and $\ell$ be the entries in the row of $\tilde G$:
\begin{equation}
 \label{eq:nikolov-table}
\begin{array}{c|c|c|c}
\text{type}&\tilde G&q_0&\ell\\ \hline
A_N(q)&\SL_{N+1}(\F_q)&q&N+1\\
B_N(q)&\Omega_{2N+1}(\F_q)&q&N\\
C_N(q)&\Sp_{2N}(\F_q)&q&N\\
D_N(q)&\Omega_{2N}^+(\F_q)&q&N\\
{}^2D_N(q)&\Omega_{2N}^-(\F_q)&q&N-1\\
{}^2A_N(q)&\SU_{N+1}(\F_{q^2})&q^2&\lceil N/2\rceil
\end{array}
\end{equation}
In untwisted type~$A$, $\tilde G=\SL_\ell(\F_{q_0})$.  In the other types, by
\cite[Sections~1 and~2]{Nikolov}, where $\ell$ and $\F_{q_0}$ are denoted $n'$
and $F_0$, the group $G$ contains a subgroup isomorphic to
$\SL_\ell(\F_{q_0})/K$ for some central subgroup $K\le\SL_\ell(\F_{q_0})$.
For a subgroup or a multiset $X$ in a group and an element $g$ of that group,
denote $X^g\defeq g^{-1}Xg$.
\begin{theorem}[Nikolov]
\label{thm:nikolov}
\begingroup
\renewcommand{\theenumi}{(\alph{enumi})}
\renewcommand{\labelenumi}{\theenumi}
\begin{enumerate}
\item\label{it:nikolov-product}
 By \cite[Theorem~1 and Section~2]{Nikolov}, if $\ell\ge2$ and $G$ is not of
 untwisted type $A$, then there are a central subgroup
 $K_0\le\SL_\ell(\F_{q_0})$, a subgroup $H_0\le G$ isomorphic to
 $\SL_\ell(\F_{q_0})/K_0$, and $g_1,\ldots,g_{200}\in G$ such that
 $G=H_0^{g_1}\cdots H_0^{g_{200}}$.
\item\label{it:nikolov-lowering}
 By \cite[Section~4, Lemma~2]{Nikolov}, for every perfect field $F$ and every
 $a\ge4$, the group $\SL_a(F)$ is a product of four conjugates of a subgroup
 isomorphic to $\SL_{a-1}(F)$.
\end{enumerate}
\endgroup
\end{theorem}

\begin{lemma}[A large type-$A$ subgroup]
\label{lem:large-A}
Let $G$ be a finite simple classical group over $\F_q$, of absolute Lie
rank $N\ge33$, and let $q_0\in\{q,q^2\}$ be as in the table
\eqref{eq:nikolov-table}.  There are integers $n\ge1$ and $c\le2^{13}$, a
subgroup $H\le G$ that is a quotient of $\SL_{3n}(\F_{q_0})$, and elements
$g_1,\ldots,g_c\in G$ such that $G=H^{g_1}\cdots H^{g_c}$ and
$|H|\ge|G|^{1/3}$.
\end{lemma}

\begin{proof}
In untwisted type $A$, take $H_0\defeq G$ and
$K_0\defeq Z(\SL_\ell(\F_{q_0}))$, so that $H_0=\SL_\ell(\F_{q_0})/K_0$ by
\eqref{eq:nikolov-table}.  In the other types, take $H_0$ and $K_0$ from
Theorem~\ref{thm:nikolov}\ref{it:nikolov-product}, and identify $H_0$ with
$\SL_\ell(\F_{q_0})/K_0$; then $G$ is a product of $200$ conjugates of $H_0$.
Let $a$ be the element of $\{\ell,\ell-1,\ell-2\}$ divisible by $3$, and denote
$n\defeq a/3$.  By Theorem~\ref{thm:nikolov}\ref{it:nikolov-lowering},
applied $\ell-a$ times (here $\ell\ge17$ by \eqref{eq:nikolov-table}), the group
$\SL_\ell(\F_{q_0})$ is a product of at most $4^2$ conjugates of a subgroup $L$
isomorphic to $\SL_a(\F_{q_0})$.  Hence $H_0$ is a product of at most $4^2$
conjugates of the image $H\defeq LK_0/K_0$ of $L$, and
$H\cong L/(L\cap K_0)$ with $L\cap K_0\le Z(L)$.  Therefore $G$ is a product
of $c\le200\cdot4^2=3200<2^{13}$ conjugates of $H$.  As $a\ge\ell-2$, by
\eqref{eq:nikolov-table},
\[
 a\ge N-1\ \text{in type }A_N,
 \quad
 a\ge N-3\ \text{in types }B_N,C_N,D_N,{}^2D_N,
 \quad
 a\ge\lceil N/2\rceil-2\ \text{in type }{}^2A_N.
\]

By the order formulas in \cite[Section~24.1 and Table~24.1]{MalleTesterman},
$|\SL_a(\F_{q_0})|=q_0^{a^2-1}\prod_{i=2}^a(1-q_0^{-i})>\tfrac12q_0^{a^2-1}$,
because $\sum_{i\ge2}2^{-i}=\tfrac12$.  Moreover,
$|L\cap K_0|\le|Z(L)|\le a$, so
\[
 \log_q|H|>(a^2-1)\log_qq_0-\log_q(2a).
\]
Here $a\le N+1$ and $q\ge2$, so $\log_q(2a)\le\log_2(2N+2)\le N$ for
$N\ge33$.  Hence, by the three bounds on $a$ and the same order formulas,
\[
\begin{array}{c|c|c|c}
\text{type}&\log_q|H|\ge&\log_q|G|\le&3\log_q|H|-\log_q|G|\ge\\ \hline
A_N(q)&N^2-3N&N^2+2N&2N^2-11N\\
B_N(q),C_N(q),D_N(q),{}^2D_N(q)&N^2-7N+8&2N^2+N&N^2-22N+24\\
{}^2A_N(q)&N^2/2-5N+6&N^2+2N+1&N^2/2-17N+17.
\end{array}
\]
All three entries of the fourth column are increasing for $N\ge33$, and at
$N=33$ they equal $1815$, $387$ and $1/2$, so all three are positive for
$N\ge33$, proving the lemma.
\end{proof}

\begin{proposition}
\label{prop:high-rank-range}
Let $G$ be a finite simple classical group of absolute Lie rank $N\ge33$.  If
$1\le m\le|G|$, then $G^m$ has a symmetric generating multiset $S$ with
$|S|\le2^{22}$ and $\gap_{G^m}(S)\ge2^{-58}$.
\end{proposition}

\begin{proof}
Take $n,c,H,g_1,\ldots,g_c$ from Lemma~\ref{lem:large-A}.  By the assumption
and Lemma~\ref{lem:large-A},
$m\le|G|\le|H|^3\le|\SL_{3n}(\F_{q_0})|^{(29-2)/9}$.  Hence, by
Proposition~\ref{prop:type-A} with $r=29$, the group
$\SL_{3n}(\F_{q_0})^m$ is a quotient of $\EL_3(R)$, and so is its quotient
$H^m$.  By Corollary~\ref{cor:standard-gap}, the image $T$ of $Y$ in $H^m$
satisfies
\[
 |T|=12\cdot30=360,
 \qquad
 \gap_{H^m}(T)\ge2^{-21}\cdot30^{-2}.
\]
By Observation~\ref{obs:direct-power-product} with $G_i=H^{g_i}$, and as
$(H^{g_i})^m=(H^m)^{(g_i,\ldots,g_i)}$,
\begin{equation}
 \label{eq:conjugate-product}
 G^m=(H^m)^{(g_1,\ldots,g_1)}\cdots(H^m)^{(g_c,\ldots,g_c)}.
\end{equation}
Apply Lemma~\ref{lem:product-transfer} to \eqref{eq:conjugate-product}, with
$T_i\defeq T^{(g_i,\ldots,g_i)}$ and $s=|T|=360$; each $T_i$ is symmetric,
so $\biguplus_{i=1}^cT_i$ is symmetric.  Since $c\le2^{13}$, the multiset
$\biguplus_{i=1}^cT_i$ has size $360c\le2^{22}$, and by
\eqref{eq:transfer-union} its gap is at least
$(2^{-21}\cdot30^{-2})/(2c^2)\ge(2^{-21}\cdot2^{-10})/(2\cdot2^{26})=2^{-58}$.
\end{proof}

\begin{remark}
\label{rem:general-r}
Let $r\ge3$.  By the proof of Proposition~\ref{prop:high-rank-range} with
$r$ in place of $29$, for every $m\le|G|^{(r-2)/27}$ the group $G^m$ has a
symmetric generating multiset $S$ with $|S|\le2^{17}(r+1)$ and
$\gap_{G^m}(S)\ge2^{-48}(r+1)^{-2}$.  Hence in high rank
Theorem~\ref{thm:amplification} is not needed for a gap of order $d^{-2}$:
if $G^m$ is $d$-generated,
then $m\le\mu_d(G)\le|G|^{d-1}$ by Observation~\ref{obs:multiplicity}, and,
for $r=27(d-1)+2$, the multiset $S$ satisfies
\[
 |S|\le2^{17}\cdot27d\le2^{22}d,
 \qquad
 \gap_{G^m}(S)\ge2^{-48}(27d)^{-2}\ge2^{-58}d^{-2}.
\]
We use $r=29$ and
Theorem~\ref{thm:amplification} instead, because with them the exponent of
$d$ in Theorem~\ref{thm:proved-global} is $\log_23$ instead of $2$.
\end{remark}

\begin{remark}
\label{rem:fixed-type}
For an untwisted type of rank at least two other than $B_2$, $C_N$ and $G_2$, one
can replace Theorem~\ref{thm:nikolov} by property~(T) of the elementary
Chevalley group of that type over a polynomial ring, proved in \cite{EJK}.
This gives a weaker form of Theorem~\ref{thm:proved-global}, with exponent
$2$ in place of $\log_23$ and with constants that depend on the rank; it is
the second row of the \hyperref[tab:known-results]{table} in the
introduction.

For a reduced irreducible root
system $\Phi$ of rank at least two and a commutative ring $A$, denote by
$E_\Phi(A)$ the group generated by the root subgroups of the simply
connected Chevalley group of type $\Phi$ over $A$.  By
\cite[Theorem~7.12]{EJK} (there, $d=r+1$), for $\Phi$ not of type $B_2$, $C_N$
or $G_2$ there is a constant $C_\Phi\ge1$ with the following property.  For
every $r\ge1$, the Steinberg group $\mathrm{St}_\Phi(\Z[x_0,\ldots,x_r])$ has
a generating set of at most $C_\Phi r$ elements with Kazhdan constant at
least $C_\Phi^{-1}r^{-1/2}$.  By \cite[Section~1.1]{EJK},
$\mathrm{St}_\Phi(\Z[x_0,\ldots,x_r])$ surjects onto
$E_\Phi(\Z[x_0,\ldots,x_r])$.

Let $q$ be a power of a prime $p$, let $G$ be the finite simple quotient of
$E_\Phi(\F_q)$, and let $G^m$ be $d$-generated.  Then $m\le\mu_d(G)\le|G|^{d-1}$ by
Observation~\ref{obs:multiplicity}, and, by the order formulas in
\cite[Section~24.1 and Table~24.1]{MalleTesterman}, $|G|<q^{2|\Phi|}$, so
$m\le q^r$ for $r\defeq2(d-1)|\Phi|$.  Let $\gamma\in\F_q$ satisfy $\F_p[\gamma]=\F_q$
and let $P_1,\ldots,P_m\in\F_q^r$ be distinct.  The ring homomorphism
$\Z[x_0,\ldots,x_r]\to\F_q^m$ that sends
$x_0$ to $(\gamma,\ldots,\gamma)$ and $x_i$ to
$\bigl((P_1)_i,\ldots,(P_m)_i\bigr)$ for $1\le i\le r$ is surjective.
Hence, as in the proof of Proposition~\ref{prop:type-A}, it induces a
surjection
$E_\Phi(\Z[x_0,\ldots,x_r])\twoheadrightarrow E_\Phi(\F_q)^m$.
The image in $G^m$ of the generating set from \cite[Theorem~7.12]{EJK},
together with the inverses, is a symmetric multiset of $s\le2C_\Phi
r\le4C_\Phi|\Phi|d$ elements, and as in Corollary~\ref{cor:standard-gap} its
gap is at least
\[
 \frac{2C_\Phi^{-2}r^{-1}}{2s}
 \ge\frac1{2C_\Phi^3r^2}
 \ge\frac1{8C_\Phi^3|\Phi|^2d^2}.
\]
For each rank there are finitely many types $\Phi$ other than $B_2$,
$C_N$ and $G_2$, so the constants may be taken to depend only on the rank.
\end{remark}

\section{Bounded rank}
\label{sec:bounded-rank}

As in Section~\ref{sec:high-rank}, we first treat type $A$, here
$\PSL_2(\F_q)$, and then all bounded rank groups.

\subsection*{Cartwright--Steger groups and \texorpdfstring{$\SL_2$}{SL2}
quotients}

Let $q$ be a prime power, and, for formal variables $y$ and $z$,
define the two algebras
\[
 R_q\defeq\F_q[y,y^{-1},(1+y)^{-1}],
 \qquad
 A_q\defeq
 \frac{\F_{q^2}[y,y^{-1},(1+y)^{-1}]\angles z}
      {\bigl(z^2-(1+y),\ z\alpha-\alpha^qz:\alpha\in\F_{q^2}\bigr)},
\]
where $\angles z$ stands for a variable that commutes with $y$ but not with
$\F_{q^2}$, and the quotient is by the two-sided ideal generated by the
listed elements.  By \cite[around equation~(4)]{LSV}, the
algebra $A_q$ is a free $R_q$-module of rank $4$.  By
\cite[Proposition~3.1 and equation~(5)]{LSV}, it splits over $\F_q((y))$,
that is, $A_q\otimes_{R_q}\F_q((y))\cong M_2(\F_q((y)))$; following
\cite[the paragraph after equation~(6), and Figure~2]{LSV}, let
$\nrd\colon A_q\to R_q$ denote the reduced norm, that is, the determinant
under that isomorphism.  By \cite[the paragraph opening Section~4]{LSV},
this splitting embeds $A_q^\times$ in $\GL_2(\F_q((y)))$ and
$A_q^\times/R_q^\times$ in $\PGL_2(\F_q((y)))$.
As $z^2=1+y$ is invertible, so is $z$; denote $b\defeq1-z^{-1}\in A_q$.
By \cite[Proposition~4.1]{LSV}, $\nrd(b)=y/(1+y)\in R_q^\times$, so $b$ is
invertible.\nolinebreak\ Denote
\[
 B_q\defeq\{ubu^{-1}R_q^\times:u\in\F_{q^2}^\times\}\subseteq A_q^\times/R_q^\times.
\]
The element $ubu^{-1}$ depends only on $u\F_q^\times$, and, by
\cite[Proposition~4.3]{LSV}, the class $ubu^{-1}R_q^\times$ determines
$u\F_q^\times$; hence $|B_q|=q+1$.  Following \cite[Definition~4.4]{LSV},
the group $\Gamma_q\defeq\angles{B_q}\le A_q^\times/R_q^\times$ is called the
\textbf{Cartwright--Steger group}.  For an ideal $J\lhd R_q$, denote by
$\varphi_J\colon A_q^\times/R_q^\times\to(A_q/JA_q)^\times/(R_q/J)^\times$ the
homomorphism induced by the quotient map $A_q\to A_q/JA_q$.  By
Wedderburn's theorem, as shown in \cite[the paragraph before
Theorem~6.3]{LSV}, $A_q/JA_q\cong M_2(R_q/J)$ when $J$ is a maximal ideal.
If moreover $R_q/J\cong\F_q$, then, by the Skolem--Noether
theorem~\cite[Section~12.6]{Pierce}, this isomorphism can be chosen to map every
$u\in\F_{q^2}\subseteq A_q$ to the matrix of multiplication by $u$ on
$\F_{q^2}\cong\F_q^2$, and we choose it so.  We call a product
$J=J_1\cdots J_k$ of distinct maximal ideals a \textbf{square-free} ideal.  For
such $J$, by the Chinese remainder theorem, the product of the isomorphisms
$A_q/J_iA_q\cong M_2(R_q/J_i)$ is an isomorphism $A_q/JA_q\cong M_2(R_q/J)$,
and the target of $\varphi_J$ is $\PGL_2(R_q/J)$.

\begin{theorem}[Lubotzky--Samuels--Vishne]
\label{thm:LSV}
Let $q$ be a prime power.
\begingroup
\renewcommand{\theenumi}{(\alph{enumi})}
\renewcommand{\labelenumi}{\theenumi}
\begin{enumerate}
\item\label{it:LSV-free}
 By \cite[Corollary~5.4]{LSV}, $B_q$ is symmetric, and $\Gamma_q$ is free of
 rank $(q+1)/2$ for $q$ odd and a free product of $q+1$ copies of $\Z/2\Z$
 for $q$ even.
\item\label{it:LSV-ramanujan}
 By \cite[Theorem~6.2]{LSV}, for every nonzero ideal $J\lhd R_q$, the graph
 $\Cay(\varphi_J(\Gamma_q);\varphi_J(B_q))$ is Ramanujan, that is,
 every eigenvalue other than $\pm(q+1)$ of its adjacency operator has
 absolute value at most $2\sqrt q$.
\end{enumerate}
\endgroup
\end{theorem}

The proof of \cite[Theorem~6.2]{LSV}\footnote{In \cite{LSV}, with their $d$
equal to $2$, our $R_q$, $A_q$,
$A_q^\times/R_q^\times$, $ubu^{-1}$, $\Gamma_q$ and $\varphi_J(\Gamma_q)$ correspond to
$R$, $\mathcal A(R)$, $G'(R)$, $b_u$, $\Gamma$ and $\Gamma/\Gamma(J)$.
Theorem~\ref{thm:LSV}\ref{it:LSV-ramanujan} is the case $d=2$ of
\cite[Theorem~6.2]{LSV}, in which the spectrum $\mathfrak S_2$ defined in
\cite[the paragraph after Remark~6.1]{LSV} is the interval
$[-2\sqrt q,2\sqrt q]$.} assumed the global Jacquet--Langlands
correspondence for function fields, which Badulescu and Roche~\cite{BR} later
proved.\footnote{For $d=2$, the Jacquet--Langlands correspondence was proved
before \cite{BR}; see \cite[Remark~3.5]{LubotzkyLie}.}

The analogue of the next proposition for a power of a prime ideal is
\cite[Theorem~6.6 and Proposition~6.7]{LSV}; for square-free ideals we
deduce it from \cite[Theorem~6.3 and Corollary~6.5]{LSV} and from the proof
of \cite[Theorem~6.3]{LSV}.

\begin{proposition}
\label{prop:LSV-square-free}
Let $q$ be a prime power and let $J=J_1\cdots J_k$ be a product of distinct
maximal ideals of $R_q$ with $|R_q/J_i|\ge4$ for every $i$.  Then
$\varphi_J(\Gamma_q)$ contains $\PSL_2(R_q/J)$, of index the order of the
class of $y/(1+y)$ in $(R_q/J)^\times/\bigl((R_q/J)^\times\bigr)^2$.
\end{proposition}

\begin{proof}
Denote by $\Lambda_q^1\le A_q^\times/R_q^\times$ the image of
$\{a\in A_q^\times:\nrd(a)=1\}$, and denote
$\Gamma_q^1\defeq\Gamma_q\cap\Lambda_q^1$, a normal subgroup of $\Lambda_q^1$
because, by \cite[Proposition~4.9]{LSV}, $\Gamma_q$ is normal in
$A_q^\times/R_q^\times$.
Each $J_i$ is a prime ideal, so, by \cite[Theorem~6.3]{LSV},
$\varphi_{J_i}(\Lambda_q^1)=\PSL_2(R_q/J_i)$, and, by
\cite[Corollary~6.5]{LSV},
$\Lambda_q^1=\Gamma_q^1\cdot\ker(\varphi_{J_i}|_{\Lambda_q^1})$.  These two
equalities together imply
\begin{equation}
 \label{eq:LSV-onto}
 \varphi_{J_i}(\Gamma_q^1)=\PSL_2(R_q/J_i).
\end{equation}
The proof of \cite[Theorem~6.3]{LSV} uses the assumption that the ideal is a
power of a prime ideal only to identify $R_q/J_i$ with the residue field of
one completion of $\F_q(y)$.
Before that, it uses the strong approximation theorem, by which
$\Lambda_q^1$ is dense in the product, over all places $\nu$ of $\F_q(y)$
other than the three places of $y$, $1/y$ and $1+y$, of the groups defined
like $\Lambda_q^1$ with the ring of integers of the completion at $\nu$ in
place of $R_q$.  Hence the same proof applies
to the $k$ places dividing $J$ at once, and proves
\begin{equation}
 \label{eq:LSV-onto-J}
 \varphi_J(\Lambda_q^1)=\PSL_2(R_q/J).
\end{equation}
By the Chinese remainder theorem, $\PSL_2(R_q/J)=\prod_i\PSL_2(R_q/J_i)$, and
$\varphi_{J_i}$ is the $i$th component of $\varphi_J$.  Hence, by
\eqref{eq:LSV-onto} and \eqref{eq:LSV-onto-J}, $\varphi_J(\Gamma_q^1)$ is a
normal subgroup of $\PSL_2(R_q/J)$ that projects onto every factor.  Each
$\PSL_2(R_q/J_i)$ is a nonabelian finite simple group because
$|R_q/J_i|\ge4$, hence a minimal normal subgroup of $\PSL_2(R_q/J)$.  By
\cite[Theorem~4.3A(i)]{DixonMortimer}, each factor either is contained in
$\varphi_J(\Gamma_q^1)$ or centralizes it.  A factor $\PSL_2(R_q/J_i)$ does
not centralize $\varphi_J(\Gamma_q^1)$, as $\varphi_J(\Gamma_q^1)$ projects onto that
factor, whose center is
trivial.\footnote{Alternatively, each factor $K_i\defeq\PSL_2(R_q/J_i)$ is
perfect, and $\varphi_J(\Gamma_q^1)$ is normal and projects onto $K_i$, so
$\varphi_J(\Gamma_q^1)\supseteq[K_i,\varphi_J(\Gamma_q^1)]=[K_i,K_i]=K_i$.}
Hence $\varphi_J(\Gamma_q)\supseteq\varphi_J(\Gamma_q^1)=\PSL_2(R_q/J)$.

We now prove the claim about the index.  Under the isomorphism
$A_q/JA_q\cong M_2(R_q/J)$, the reduced norm modulo $J$ is the
determinant, which induces a homomorphism from $\PGL_2(R_q/J)$ onto
$(R_q/J)^\times/\bigl((R_q/J)^\times\bigr)^2$ with kernel $\PSL_2(R_q/J)$.  Now
$\Gamma_q=\angles{B_q}$, and $\nrd(ubu^{-1})=\nrd(b)=y/(1+y)$ for every $u$,
because the reduced norm is invariant under conjugation.  Therefore the
image of
$\varphi_J(\Gamma_q)$ in $(R_q/J)^\times/\bigl((R_q/J)^\times\bigr)^2$ is the subgroup
generated by the class of $y/(1+y)$, and
${[\varphi_J(\Gamma_q):\PSL_2(R_q/J)]}$ is the order of that class.
\end{proof}

For a prime power $q$, denote $m_q\defeq(q-3)/2$ for $q$ odd and
$m_q\defeq q-2$ for $q$ even.\footnote{The fields of order $2$ and $3$ are
excluded from the next proposition and are treated, together with $\F_4$ and
$\F_5$, in Corollary~\ref{cor:small-fields}.}

\begin{proposition}[Square-free Cartwright--Steger quotient]
\label{prop:LSV-product}
Let $q\ge4$ be a prime power.  There is a square-free ideal $J\lhd R_q$ with
$R_q/J\cong\F_q^{m_q}$ for which $\varphi_J(\Gamma_q)=\PSL_2(\F_q)^{m_q}$, and
$\varphi_J(B_q)$ is a
symmetric generating multiset of $q+1$ elements.
Moreover, for $G\defeq\PSL_2(\F_q)^{m_q}$, the graph $\Cay(G;\varphi_J(B_q))$
is Ramanujan, and every eigenvalue of $M_{\varphi_J(B_q)}^{\Reg_G}$ on
$\mathbf1^\perp$ is at most $2\sqrt q/(q+1)$.
\end{proposition}

\begin{proof}
Let $\mathcal A$ be the set of $a\in\F_q\setminus\{0,-1\}$ for which
$a/(1+a)$ is a square, and denote
$J\defeq\bigl(\prod_{a\in\mathcal A}(y-a)\bigr)\lhd R_q$.  The map $a\mapsto a/(1+a)$ is a bijection from $\F_q\setminus\{0,-1\}$ onto
$\F_q\setminus\{0,1\}$, so $|\mathcal A|=m_q$: if $q$ is odd, exactly
$(q-3)/2$ elements of $\F_q\setminus\{0,1\}$ are squares, and if $q$ is even,
all $q-2$ are.  For $a\in\mathcal A$, since $a\ne0,-1$, the ideal $(y-a)$ is
maximal in $R_q$, with $R_q/(y-a)\cong\F_q$; these ideals are distinct, so
$J$ is a product of $m_q$ distinct maximal ideals, and, by the Chinese
remainder theorem, $R_q/J\cong\F_q^{m_q}$.  The image of $y/(1+y)$ in
$R_q/(y-a)$ is $a/(1+a)$, which is a square by the choice of $\mathcal A$, so
the class of $y/(1+y)$ in $(R_q/J)^\times/\bigl((R_q/J)^\times\bigr)^2$ is
trivial, and, by Proposition~\ref{prop:LSV-square-free}, which applies as
$q\ge4$,
\[
 \varphi_J(\Gamma_q)=\PSL_2(R_q/J)=\PSL_2(\F_q)^{m_q}.
\]
By Theorem~\ref{thm:LSV}\ref{it:LSV-free}, the set $B_q$ is symmetric, and
it generates $\Gamma_q$, so $\varphi_J(B_q)$ is a symmetric generating
multiset of $q+1$ elements.  As $(q+1)M_{\varphi_J(B_q)}^{\Reg_G}$ is the
adjacency operator of $\Cay(G;\varphi_J(B_q))$, by
Theorem~\ref{thm:LSV}\ref{it:LSV-ramanujan} it remains to exclude the
eigenvalue $q+1$ of this operator on $\mathbf1^\perp$, and this holds as the
graph is connected.
\end{proof}

\begin{theorem}[{Lubotzky~\cite[Theorem~3.7]{LubotzkyLie}}]
\label{thm:lubotzky}
The family $\{\PSL_2(\F_q)\}$, where $q$ runs over the prime powers, is an
expanding family.  More explicitly, there is an absolute constant
$\varepsilon_0>0$ such that, for every prime power $q$, there is a multiset
$S_q^+$ of four elements of $\PSL_2(\F_q)$ for which
$S_q\defeq S_q^+\uplus(S_q^+)^{-1}$ satisfies
$\gap_{\PSL_2(\F_q)}(S_q)\ge\varepsilon_0$.
\end{theorem}

By the discrete Cheeger inequality \cite[Proposition~4.2.4]{LubotzkyBook},
Theorem~\ref{thm:lubotzky} follows from \cite[Theorem~3.7]{LubotzkyLie},
which proves a combinatorial version of expansion.  The proof of
Proposition~\ref{prop:bounded-rank-one} is the proof of
\cite[Theorem~3.7]{LubotzkyLie}, adapted to the power $\PSL_2(\F_q)^{m_q}$.

\begin{proposition}
\label{prop:bounded-rank-one}
The families $\{\PSL_2(\F_q)^{m_q}\}$ and $\{\SL_2(\F_q)^{m_q}\}$, where
$q\ge4$ runs over the prime powers, are expanding families.  More
explicitly, there is an absolute constant $\varepsilon_1>0$ with the
following property.  Let $q\ge4$ be a prime power, and let $G$ be
$\PSL_2(\F_q)^{m_q}$ or $\SL_2(\F_q)^{m_q}$.  Then $G$ has a symmetric
multiset $U$ with $|U|\le12$ and $\gap_G(U)\ge\varepsilon_1$.
\end{proposition}

\begin{proof}
Consider first $G=\PSL_2(\F_q)^{m_q}$.  By Theorem~\ref{thm:lubotzky},
there are an absolute constant $\varepsilon_0>0$ and a multiset $S_q^+$ of
four elements of $\PSL_2(\F_q)$ such that $S_q\defeq S_q^+\uplus(S_q^+)^{-1}$
satisfies $\gap_{\PSL_2(\F_q)}(S_q)\ge\varepsilon_0$.  Denote by
$\Delta\colon\PGL_2(\F_q)\to\PGL_2(\F_q)^{m_q}$ the diagonal embedding
$g\mapsto(g,\ldots,g)$, and let $K\defeq\Delta(\PSL_2(\F_q))$ and
$S_K\defeq\Delta(S_q)$; then $\gap_K(S_K)\ge\varepsilon_0$.

Let $T\le\PGL_2(\F_q)$ be the image of the map sending $u\in\F_{q^2}^\times$
to the matrix of multiplication by $u$ on $\F_{q^2}\cong\F_q^2$; it is
called a \textbf{nonsplit torus}.  Let $J$ be as in Proposition~\ref{prop:LSV-product}, so that
$R_q/J\cong\F_q^{m_q}$, and denote $S_G\defeq\varphi_J(B_q)$.  By the choice of the isomorphisms before
Theorem~\ref{thm:LSV}, the homomorphism $\varphi_J$ maps the image of
$\F_{q^2}^\times$ in $A_q^\times/R_q^\times$ onto $\Delta(T)$.  Hence, by the
definition of $B_q$, conjugation by $\Delta(T)$ permutes $S_G$
transitively.  As $\PSL_2(\F_q)$ has index at most two in $\PGL_2(\F_q)$, the
elements of $S_G$ form at most two orbits under conjugation by $\Delta(T\cap\PSL_2(\F_q))\le K$.  Let
$B_0\subseteq S_G$ consist of one element of each orbit, and denote
$S_G'\defeq B_0\uplus B_0^{-1}$.  Then, for every $x\in S_G$, there are
$t\in K$ and $x'\in S_G'$ such that $x=t^{-1}x't$.  By
Proposition~\ref{prop:LSV-product} and the definition \eqref{eq:gap-def} of
the gap, the multiset $S_G$ is symmetric and
$\gap_G(S_G)\ge1-2\sqrt q/(q+1)\ge1/5$, as $2\sqrt q/(q+1)$ decreases in
$q\ge1$ and equals $4/5$ at $q=4$.  Hence, by
Lemma~\ref{lem:conjugate-transfer} with $\varepsilon_K=\varepsilon_0$ and
$\varepsilon_G=1/5$, the multiset $U\defeq S_K\uplus S_G'$, of size
$8+2|B_0|\le12$, satisfies
\[
 \gap_G(U)
 \ge\frac{\varepsilon_G}{\bigl(16/(|S_K|\,\varepsilon_K)+2\bigr)|U|}
 \ge\frac{1/5}{\bigl(16/(8\varepsilon_0)+2\bigr)\cdot12}
 =\frac1{120\,(1+\varepsilon_0^{-1})}.
\]

The group $\SL_2(\F_q)^{m_q}$ is a perfect central extension of the perfect
group $G$.  Let $\widetilde U$ be the multiset in $\SL_2(\F_q)^{m_q}$ that
consists of one preimage, under the quotient map, of each element of
$\Delta(S_q^+)\uplus B_0$, and of the inverses of these preimages; then the
quotient map takes $\widetilde U$ onto $U$.  By
\cite[Proposition~4.2.5]{LubotzkyBook} (from the gap to combinatorial
expansion), \cite[Proposition~2.3]{LubotzkyLie} (from $G$ to its perfect
central extension) and \cite[Proposition~4.2.4]{LubotzkyBook} (back to the
gap), $\gap_{\SL_2(\F_q)^{m_q}}(\widetilde U)\ge\varepsilon'$ for an absolute
constant $\varepsilon'>0$.  Hence the proposition holds with
${\varepsilon_1\defeq\min\bigl\{(120(1+\varepsilon_0^{-1}))^{-1},\varepsilon'\bigr\}}$.
\end{proof}

\subsection*{From \texorpdfstring{$\SL_2$}{SL2} to bounded rank}

\begin{proposition}[Bounded rank]
\label{prop:bounded-rank}
There is an absolute constant $C_1>0$ with the following property.  Let $G$
be a finite simple group of Lie type over $\F_q$ with $q\ge7$, other than
${}^2B_2(q)$, of absolute Lie rank below $33$.\footnote{For a twisted group,
$q$ is the prime power in its Dynkin-type label, as in
Section~\ref{sec:preliminaries} and \cite[Table~1]{LNS}; for example,
${}^2A_N(q)=\PSU_{N+1}(\F_{q^2})$.}
Then every $d$-generated power $G^m$ has a symmetric generating multiset $S$
with $|S|\le C_1d$ and $\gap_{G^m}(S)\ge C_1^{-1}d^{-\log_23}$.
\end{proposition}

\begin{proof}
By \cite[Theorem~1.1 and Table~1]{LNS}, $G$ is a product of $c_0$ subgroups
isomorphic to $\SL_2(\F_q)$ or $\PSL_2(\F_q)$.  If $G$ is untwisted, then
$c_0=5|\Phi^+|$, where $\Phi^+$ is its set of positive roots.  For the
twisted types ${}^2A_N(q)$, ${}^2D_N(q)$, ${}^3D_4(q)$, ${}^2E_6(q)$,
${}^2F_4(q)$ and ${}^2G_2(q)$, the number $c_0$ is given in
\cite[Table~1]{LNS} or, for the few groups omitted there, in \cite[Remark~3]{LNS}.  For rank below $33$ and $q\ge7$, all these numbers
are at most $2^{13}$.  For each of these $c_0$ subgroups $H$, apply
Proposition~\ref{prop:bounded-rank-one} to $H^{m_q}$, adding copies of the
identity to get $12$ elements.  By Observation~\ref{obs:direct-power-product}
and Lemma~\ref{lem:product-transfer}, the group $G^{m_q}$ then has a
symmetric multiset of $12c_0\le12\cdot2^{13}$ elements, with gap at least an
absolute constant $\eta_1>0$.

Since the rank is bounded, by the order formulas
\cite[Section~24.1 and Table~24.1]{MalleTesterman}, there is an absolute
constant $\beta_0$ with $|G|\le q^{\beta_0}$, and
$m_q\ge q^{1/3}\ge|G|^{1/(3\beta_0)}$ as $q\ge7$.  Now apply
Corollary~\ref{cor:reach} with $(m_0,\beta,s,\eta)=(m_q,3\beta_0,12c_0,\eta_1)$.
\end{proof}

\section{Proofs of the main theorems}
\label{sec:proofs}

In this section we collect the results of the previous sections and prove
Theorems~\ref{thm:proved-global} and~\ref{thm:huge-powers}.
Theorem~\ref{thm:proved-global} also follows from
Theorem~\ref{thm:huge-powers}, by Corollary~\ref{cor:reach} with $m_0=|G|$
and $\beta=1$.  We prove Theorem~\ref{thm:proved-global} first because a
direct proof of Theorem~\ref{thm:huge-powers} would still use
Theorem~\ref{thm:amplification} in bounded rank, through
Corollary~\ref{obs:balanced-iteration}.

\begin{corollary}
\label{cor:small-fields}
Let $\mathcal F$ be a finite set of nonabelian finite simple groups.  There
is $C_2>0$ such that, for every $G\in\mathcal F$, every $d$-generated power
$G^m$ has a symmetric generating multiset $S$ with $|S|\le C_2d$ and
$\gap_{G^m}(S)\ge C_2^{-1}d^{-\log_23}$.  Indeed, $\gap_{G^2}(G^2)=1$ and
$2=|G|^{1/\log_2|G|}$, so Corollary~\ref{cor:reach} applies with
$(m_0,\beta,s,\eta)=(2,\log_2|G|,|G|^2,1)$, and $C_2$ is the maximum over
$\mathcal F$ of the resulting constants.
\end{corollary}

\begin{proof}[Proof of Theorem~\ref{thm:proved-global}]
Let $G\in\mathcal L$.  By \cite[Table~22.1]{MalleTesterman}, every
exceptional Lie type has absolute Lie rank at most $8$, so a group of Lie
type of absolute Lie rank at least $33$ is classical.  For such a $G$, by
Proposition~\ref{prop:high-rank-range}, the group $G^{|G|}$ has a symmetric
generating multiset of size at most $2^{22}$ and gap at least $2^{-58}$.
Hence, by Corollary~\ref{cor:reach} with
$(m_0,\beta,s,\eta)=(|G|,1,2^{22},2^{-58})$, every $d$-generated power $G^m$
has a symmetric generating multiset $S$ with
\[
 |S|\le2^{23}d,
 \qquad
 \gap_{G^m}(S)\ge2^{-60}d^{-\log_23}.
\]
The remaining groups of Lie type over $\F_q$ with $q\ge7$ are covered by
Proposition~\ref{prop:bounded-rank}, with its constant $C_1$.  Apply
Corollary~\ref{cor:small-fields} (with a constant $C_2$) to the finite set
of the groups in $\mathcal L$ that are sporadic or of Lie type over $\F_q$
with $q\le5$ and absolute Lie rank below $33$.  Then
$C\defeq\max\{2^{60},C_1,C_2\}$ works in all three cases.
\end{proof}

\begin{proof}[Proof of Theorem~\ref{thm:huge-powers}]
Let $G\in\mathcal L$.  By Observation~\ref{obs:multiplicity}, we have
$\mu_3(G)\ge|G|$, so $G^{|G|}$ is generated by $3$ elements, and
Theorem~\ref{thm:proved-global} with $d=3$ applies.
\end{proof}

\begin{corollary}
\label{cor:sabatini}
The family $\{\PSL_2(\F_p)^p\}_{p\text{ prime}}$ is an expanding family; this answers \cite[Question~4.5]{Sabatini}.  Indeed, for
$p\ge7$ the group $\PSL_2(\F_p)$ is in $\mathcal L$, and $\PSL_2(\F_p)^p$ is a
quotient of $\PSL_2(\F_p)^{|\PSL_2(\F_p)|}$, so Theorem~\ref{thm:huge-powers}
and \eqref{eq:quotient-gap} apply.
\end{corollary}

\begin{remark}
\label{rem:exponent}
The proof of Theorem~\ref{thm:proved-global} gives more than its statement:
for each $G\in\mathcal L$, the exponent $\log_23$ may be replaced by
\[
 \nu_G\defeq\log_2\frac{2\mindim(G)}{\mindim(G)-1}
 =1+\log_2\left(1+\frac1{\mindim(G)-1}\right),
\]
and $\nu_G\to1$ as $|G|\to\infty$ over $\mathcal L$.  We state
Theorem~\ref{thm:proved-global} with $\log_23$ because one exponent then
works for every $G\in\mathcal L$: by Fact~\ref{obs:mindim},
$\nu_G\le\log_23$, with equality when $\mindim(G)=3$.

To see that $\nu_G$ may replace $\log_23$, take $A=B$ in
Theorem~\ref{thm:amplification}: the factor
$\frac12\bigl(1-\frac1{\mindim(G)}\bigr)$ in the first inequality of
\eqref{eq:amplification-gap} equals $2^{-\nu_G}$, and in
Corollary~\ref{cor:square-simple} it is replaced by $\frac13\le2^{-\nu_G}$, as
$\mindim(G)\ge3$ by Fact~\ref{obs:mindim}.  In the proof of
Corollary~\ref{obs:balanced-iteration}, Corollary~\ref{cor:square-simple} is
applied $j$ times, and in the proof of Corollary~\ref{cor:reach}, $j$ is chosen
with $2^j<2\beta d$.  With $2^{-\nu_G}$ kept, the gap bound becomes
$\eta\,(2^j)^{-\nu_G}>\eta\,(2\beta d)^{-\nu_G}\ge\eta\,(2\beta)^{-\log_23}d^{-\nu_G}$,
as $2\beta\ge1$ and $\nu_G\le\log_23$.  All three cases in the proof of
Theorem~\ref{thm:proved-global} use Corollary~\ref{cor:reach}, so
$\gap_{G^m}(S)\ge C^{-1}d^{-\nu_G}$ with the same constant $C$.

For the limit, $\nu_G\le1+\frac1{(\mindim(G)-1)\ln2}$, because
$\log_2(1+x)\le x/\ln2$ for $x\ge0$.  A nontrivial irreducible representation
of $G$ is faithful, so $G$ embeds in $\GL_{\mindim(G)}(\C)$, and, by Jordan's
theorem, in the form proved in \cite{Collins}, $G$ has an abelian normal subgroup of index bounded
in terms of $\mindim(G)$.  That subgroup is trivial, so $|G|$ is bounded in
terms of $\mindim(G)$, and hence $\nu_G\to1$ as $|G|\to\infty$ over
$\mathcal L$.
\end{remark}

In the following remark we adapt to powers the second proof of
\cite[Corollary~4.4]{LW}, a proof that Lubotzky and Weiss credited to Yael
Luz.

\begin{remark}[Upper bounds]
\label{rem:upper-bound}
For $G_N\defeq\PSL_{N+1}(\F_2)$ and $d\ge2$, no choice of $d$ generators of
each group $G_N^{\mu_d(G_N)}$, for $N\ge3$, makes these groups an expanding
family.  More precisely, let $N\ge3$, denote $G\defeq G_N$, and let
$y_1,\ldots,y_d$ generate $G^{\mu_d(G)}$; then
$\Cay(G^{\mu_d(G)};\{y_1,\ldots,y_d\})$ is not a $2/(d(N+1))$-expander.
Indeed, denote by $e_1,\ldots,e_{N+1}$ the standard basis of $\F_2^{N+1}$,
by $\mathbb P\defeq\mathbb P(\F_2^{N+1})$ the projective space, by
$[x]\in\mathbb P$ the point of $0\ne x\in\F_2^{N+1}$, and, for
$X\subseteq\mathbb P$, by $\mathbf1_X\in\ell^2(\mathbb P)$ the indicator
function of $X$.  Denote $\mathcal E\defeq\{[e_1],\ldots,[e_{N+1}]\}$, let
$a\defeq I_{N+1}+E_{12}$, and let $b$ be the permutation matrix with
$be_k=e_{k+1}$ for $k\le N$ and $be_{N+1}=e_1$; by \cite[Section~3]{Taylor},
$a$ and $b$ generate $G$.  Then
$a^{-1}=a$, $b\mathcal E=\mathcal E$, and the symmetric difference of
$\mathcal E$ and $a\mathcal E$ is $\{[e_2],[e_1+e_2]\}$.  Denote by $S$ the
multiset $\{y_1^{\pm1},\ldots,y_d^{\pm1}\}$, of $2d$ elements.  By
Observation~\ref{obs:unique-tuple}, there are an index
$i\le\mu_d(G)$ and an $\alpha\in\Aut(G)$ such that the $i$th coordinates of
$y_1,\ldots,y_d$ are $\alpha(a),\alpha(b),1,\ldots,1$.  Let
$\pi\colon G^{\mu_d(G)}\to\GL(\ell^2(\mathbb P))$ be the permutation
representation of $G$, precomposed with
$\alpha^{-1}$ and with the projection of $G^{\mu_d(G)}$ onto coordinate $i$,
so that $y_1$ and $y_2$ act as $a$ and $b$ and $y_3,\ldots,y_d$ act
trivially.  Let
$v\defeq\mathbf1_{\mathcal E}-\frac{N+1}{|\mathbb P|}\mathbf1_{\mathbb P}$.
Both $y_1^{\pm1}$ map $\mathbf1_{\mathcal E}$ to $\mathbf1_{a\mathcal E}$, at
squared distance $\|\mathbf1_{a\mathcal E}-\mathbf1_{\mathcal E}\|^2=2$, while
$y_2^{\pm1},\ldots,y_d^{\pm1}$ fix $\mathbf1_{\mathcal E}$, and every $y_k$
fixes $\mathbf1_{\mathbb P}$.  Hence, by the formula \eqref{eq:dirichlet} for
$E_S^\pi$, we have $E^\pi_S(v)=\frac{2+2}{2\cdot2d}=\frac1d$.  The entries of
$v$ sum to $0$ and $G$ acts transitively on $\mathbb P$, so
$P_{G^{\mu_d(G)}}v=0$.  Moreover,
$\|v\|^2=(N+1)\bigl(1-(N+1)/|\mathbb P|\bigr)>(N+1)/2$, as
$|\mathbb P|=2^{N+1}-1>2(N+1)$.  Hence, by Observation~\ref{obs:gap-transfer},
$\gap_{G^{\mu_d(G)}}(S)\le E^\pi_S(v)/\|v\|^2<2/(d(N+1))$.

As for $\PSL_{N+1}(\F_2)$ above, for the alternating groups the analogue of
Question~\ref{q:maximal-power} has a negative answer, because some generating $d$-tuples of $\Alt(n)$ have
spectral gap tending to $0$ as $n\to\infty$ (see
\cite[Example~5.6]{Jez}).  Jezernik proved~\cite[Theorem~1.13]{Jez} a
stronger, robust version: if $\Alt(n)^{f(n)}$ has a generating $d$-tuple
whose Cayley graph has spectral gap at least $\varepsilon$, then
$\log f(n)\le(d-1-\varepsilon/160+o(1))\,n\log n$, while
$\log\mu_d(\Alt(n))\sim(d-1)\,n\log n$.  He also
asked~\cite[Question~1.12]{Jez} whether there are $\delta>0$ and
$f(n)\ge|\Alt(n)|^\delta$ such that the groups $\Alt(n)^{f(n)}$ have
generating tuples of bounded size whose Cayley graphs have spectral gap
bounded away from $0$.  By Corollary~\ref{cor:reach}, a positive answer is
equivalent to the analogue of Theorem~\ref{thm:proved-global} for the
alternating groups.
\end{remark}

\appendix
\section*{Glossary of Notation}
\addcontentsline{toc}{section}{Glossary of Notation}

\begingroup
\setlength{\parskip}{0pt}
\footnotesize
\begin{longtable}{@{}l@{\quad}p{0.62\textwidth}@{}}
\multicolumn{2}{@{}l}{\textit{Sets, groups and maps}}\\[2pt]
$\F_q$ & the field with $q$ elements, $q=p^f$ with $p$ prime\\
$|X|$ & cardinality, counting multiplicity in a multiset\\
$A\uplus B$ & the multiset in which multiplicities are summed\\
$X^g$ & $g^{-1}Xg$, for a subgroup or a multiset $X$\\
$D_d(G)$ & the set of generating $d$-tuples of $G$\\
$\mu_d(G)$ & $|D_d(G)|/|\Aut(G)|$, the number of factors of the largest $d$-generated power, \eqref{eq:hall-count}\\
$\mathcal L$ & the nonabelian finite simple groups other than the alternating and the Suzuki groups\\
$A_N(q),\ N$ & Dynkin-type labels, and the absolute Lie rank, Section~\ref{sec:preliminaries}\\
$\Cay(G;S)$ & Cayley graph\\
$\gap_G(S)$ & normalized one-sided spectral gap, \eqref{eq:gap-def}\\
$\row,\col$ & the embeddings $G^m\to G^{mn}$ and $G^n\to G^{mn}$ along the rows and the columns of an $m$-by-$n$ table\\
$G_{\row},G_{\col}$ & their images $\row(G^m)$ and $\col(G^n)$\\
$\Sigma_{A,B}$ & $\row(A)\uplus\col(B)$, a multiset in $G^{mn}$\\[6pt]

\multicolumn{2}{@{}l}{\textit{Representations}}\\[2pt]
$V,\pi$ & a finite-dimensional complex Hilbert space and a unitary representation on it\\
$\angles{\cdot,\cdot},\ \|\cdot\|$ & inner product and its norm; $\angles{X}$ is also the subgroup generated by $X$, and $\Z\angles{x_1,\ldots,x_r}$ a free ring\\
$\opnorm{\cdot}$ & operator norm\\
$I$ & the identity operator; $I_n$ is the identity matrix\\
$\ell^2(G),\ \Reg_G$ & the left regular representation of $G$\\
$\mathbf 1$ & the constant functions in $\ell^2(G)$\\
$M_S^\pi$ & averaging operator of $S$\\
$E_S^\pi(v)$ & the quadratic form $\angles{(I-M_S^\pi)v,v}$, \eqref{eq:dirichlet}\\
$V^G,\ P_G$ & the $G$-fixed subspace and the orthogonal projection onto it\\
$\mindim(G)$ & the minimal dimension of a nontrivial irreducible representation of $G$\\[6pt]

\multicolumn{2}{@{}l}{\textit{High rank (Section~\ref{sec:high-rank})}}\\*[2pt]
$R,\ r$ & the free associative ring $\Z\angles{x_1,\ldots,x_r}$, and its number of variables\\
$M_n(\F_q)$ & the algebra of $n$-by-$n$ matrices\\
$E_{ij},\ e_{ij}(a)$ & the matrix with a single $1$, in position $(i,j)$, and $I_k+aE_{ij}$\\
$\EL_k(A)$ & the group generated by the $e_{ij}(a)$ with $a\in A$\\
$Y$ & a generating set of $\EL_3(R)$, of $12(r+1)$ elements\\
$\kappa_r$ & a lower bound for its Kazhdan constant\\
$H$ & a large type-$A$ subgroup of $G$, a quotient of $\SL_{3n}(\F_{q_0})$, Lemma~\ref{lem:large-A}\\
$q_0,\ell$ & the field size and the degree of $\SL_\ell(\F_{q_0})$ in Theorem~\ref{thm:nikolov}\ref{it:nikolov-product}, \eqref{eq:nikolov-table}\\[6pt]

\multicolumn{2}{@{}l}{\textit{Bounded rank (Section~\ref{sec:bounded-rank})}}\\[2pt]
$R_q,\ y$ & $\F_q[y,y^{-1},(1+y)^{-1}]$\\
$A_q,\ z$ & $\F_{q^2}[y,y^{-1},(1+y)^{-1}]\angles z\big/\bigl(z^2-(1+y),\ z\alpha-\alpha^qz:\alpha\in\F_{q^2}\bigr)$\\
$\nrd$ & reduced norm $A_q\to R_q$\\
$b,\ B_q$ & $1-z^{-1}$, and the images in $A_q^\times/R_q^\times$ of the $q+1$ conjugates $ubu^{-1}$\\
$\Gamma_q$ & the Cartwright--Steger group $\angles{B_q}$\\
$J,J_i$ & an ideal $J\lhd R_q$, and the distinct maximal ideals with $J=J_1\cdots J_k$\\
$\varphi_J$ & reduction of coefficients modulo $J$\\
$m_q$ & the number of factors of the congruence quotient, Proposition~\ref{prop:LSV-product}\\
\end{longtable}
\endgroup

\end{document}